\documentclass[pdflatex,sn-mathphys-num]{sn-jnl}

\usepackage{graphicx}%
\usepackage{multirow}%
\usepackage{amsmath,amssymb,amsfonts}%
\usepackage{amsthm}%
\usepackage{mathrsfs}%
\usepackage[title]{appendix}%
\usepackage{xcolor}%
\usepackage{textcomp}%
\usepackage{manyfoot}%
\usepackage{booktabs}%
\usepackage{algorithm}%
\usepackage{algorithmicx}%
\usepackage{algpseudocode}%
\usepackage{listings}%

\usepackage{physics}

\usepackage{pifont}

\usepackage{tikz}
\usepackage{xcolor}
\usetikzlibrary{patterns}

\theoremstyle{thmstyleone}%
\theoremstyle{thmstyletwo}%

\theoremstyle{thmstylethree}%

\begin{document}

\title[Article Title]{A semi-implicit double-point Material Point Method for both free-surface flow and seepage in deformable porous media}


\author*[1,2]{\fnm{Mian} \sur{Xie}}\email{mian.xie@lsbu.ac.uk}

\author[3]{\fnm{Pedro} \sur{Navas}}\email{pedro.navas@upm.es}

\author[2]{\fnm{Susana} \sur{L\'opez-Querol}}\email{s.lopez-querol@ucl.ac.uk}

\affil*[1]{\orgdiv{School of Engineering and Design}, \orgname{London South Bank University}, \orgaddress{\street{103 Borough Road}, \city{London}, \postcode{SE1 0AA}, \country{UK}}}

\affil*[2]{\orgdiv{Department of Civil, Environmental and Geomatic Engineering}, \orgname{University College London}, \orgaddress{\street{Gower Street}, \city{London}, \postcode{WC1E 6BT}, \country{UK}}}

\affil[3]{\orgdiv{ETSI Caminos Canales y Puertos}, \orgname{Universidad Polit\'ecnica de Madrid}, \orgaddress{\street{Calle Profesor Aranguren 3}, \city{Madrid}, \postcode{28040}, \country{Spain}}}

\abstract{A new semi-implicit, two-phase, double-point formulation of the Material Point Method (MPM) for soil–water interaction with seepage and free-surface flows under large deformation is presented in this paper. The approach advances the water phase implicitly while keeping the soil phase explicit, enabling stable, efficient time integration in problems that involve rapid seepage and strong free-surface motion. The proposed framework models high-Reynolds-number
interphase drag through a non-linear Darcy's law implemented for the first time within an incremental fractional step MPM formulation without enlarging the implicit solve. This methodology also enhances the numerical stability for fast flows and wave breaking via a hyperelastic constitutive treatment of slightly compressible viscous water, and mitigates spurious oscillations through a new stabilisation approach for the velocity. Robustness of soil-water interface is achieved by combining nodal-based, free-surface detection, suited for higher-order spline functions with smooth porosity-permeability transitions that avoid constitutive divergence at sharp material boundaries. Validation against laboratory benchmark cases reported in the literature, including pure-water dam break, dam-break seepage through a porous barrier, two granular-collapse tsunami experiments, and a dam-break wave over a movable granular bed, shows accurate and stable free-surface evolution, pressure time histories, seepage fronts, and wave-gauge records. Using an advanced critical-state soil model (NorSand) further improves the reproduction of granular flow kinematics. The results demonstrate that the proposed formulation is a reliable and computationally efficient tool for geotechnical hazards involving intense soil–water coupling, seepage, sediment transport and free water. }




\maketitle

\section{Introduction}\label{sec1}

{Geotechnical hazards such as landslide-generated tsunamis, dam-break flows through or over embankments, and impacts of debris on infrastructures, involve intense two-way coupling between deformable saturated soils and free-surface water. Modelling these phenomena is computationally challenging because several demanding features coexist: large deformations of the solid phase, rapid free-surface evolution with wave breaking, seepage through porous media with sharp porosity gradients, and history-dependent constitutive behaviour of the geomaterial. A reliable numerical framework must capture all of these aspects simultaneously.}

{In the existing literature, a variety of computational approaches have been applied to fluid-granular interaction problems of this kind. Smoothed Particle Hydrodynamics (SPH) naturally tracks free surfaces and has been used to model dam-break flows and landslide-water interaction with elastoplastic soil models \citep{Bui2008,Bui2021,Wang2022sph}; however, tensile instability and accurate stress integration for history-dependent solids remain active challenges \citep{Bui2021}. The Moving Particle Semi-implicit (MPS) method shares similar Lagrangian advantages for free-surface tracking \citep{Koshizuka1996,Jandaghian2021}, but has typically been coupled with simplified rheologies rather than advanced constitutive models. Finite Volume Method (FVM) solvers such as OpenFOAM, combined with the Volume of Fluid technique, have been successfully applied to granular-fluid flows using the $\mu(I)$ rheology \citep{Rauter2022,Clous2023}; yet incorporating history-dependent constitutive models with evolving state variables (e.g.\ void ratio, over-consolidation) is not straightforward within a purely Eulerian framework. Coupled DEM-CFD approaches provide grain-scale resolution of fluid-particle interaction \citep{Zhao2013,Shan2014,Cheng2018} but remain computationally prohibitive at engineering scales.}

{Among these alternatives, the} Material Point Method (MPM) has become a powerful framework for simulating large deformations, contacts, and history-dependent constitutive response in geomaterials \citep{Sulsky1994, Soga2016, Augarde2021}. By combining a Lagrangian set of material points with an auxiliary background mesh, MPM avoids the distortion of the elements while retaining a clear kinematic and constitutive description at the material point level—an attractive alternative to purely Eulerian or Lagrangian finite element strategies for problems involving large strains \citep{Zhang2016}. {Crucially, MPM naturally accommodates advanced elastoplastic and critical-state constitutive models at material points while tracking the full deformation history, a combination that is difficult to achieve with the Eulerian or meshfree methods mentioned above. This makes it particularly well suited for soil-water coupled problems where both large-strain geomechanics and free-surface hydrodynamics must be captured within a single framework.}

{Within MPM, the consideration of coupling between the soil skeleton and the pore water is essential to capture consolidation, seepage, and dynamic soil-fluid interactions.} Existing two-phase MPM formulations can be broadly classified into two categories:
\emph{(i) single-point approaches}, in which a \emph{single} set of MP carries the state of both phases (soil stress and pore pressure are co-located and advected with the skeleton) \citep{Zhang2009,Alonso2010,Zabala2011,Jassim2013,Iaconeta2019,Zheng2022,Kularathna2021,Yuan2023, Pretti2025}; and
\emph{(ii) double-point approaches}, in which \emph{two} distinct MP sets represent soil and water, respectively, and interact through momentum exchange and mass compatibility \citep{Abe2014,Bandara2015,Liu2017,Yamaguchi2020,Chandra2024b,He2024}.
Historically, most soil–water coupled MPMs were developed as single-point formulations, reflecting their close methodological affinity with the Finite Element Method (FEM), which naturally colocates solid and fluid unknowns at quadrature points. In contrast, a FEM analogue of a discretisation \emph{double-point}, in which solid and fluid phases are advected by different material points, is not straightforward and does not have a direct counterpart in standard FEM practice. 

While single-point MPMs have been successful for many consolidation and deformation problems \citep{Zhang2009,Alonso2010,Zabala2011,Jassim2013,Iaconeta2019,Zheng2022,Kularathna2021,Yuan2023, Pretti2025}, they face fundamental limitations when strong phase-relative motion arises. In such cases, double-point MPMs, in which the water phase moves independently of the skeleton, are better suited to model phenomena such as free-surface flows interacting with highly deformable soils, rapid seepage with sharp porosity gradients, erosion, and landslide-generated waves \citep{Yamaguchi2020,Chandra2024b,He2024}. Indeed, double-point formulations have been shown to be not only more general, but also \emph{more accurate} even in problems that can, in principle, be treated by single-point approaches, owing to their consistent treatment of relative kinematics and interphase drag \citep{Xie2025}.

Explicit two-phase MPMs are conceptually simple, but they suffer from severe time-step restrictions controlled by the Courant–Friedrichs–Lewy (CFL) condition \citep{Courant1967} of (stiff) fluid compressibility and, under low permeability, by hydraulic diffusion constraints \citep{Mieremet2016}; this can make simulations prohibitively expensive. Semi-implicit formulations, in which the behaviour of the water phase is implicitly determined while the soil phase remains explicit, alleviate these restrictions and improve the pressure stability by satisfying a discrete inf–sup condition during the projection step \citep{Yamaguchi2020,Pan2021,Chandra2024b, He2024}. The fractional step method \citep{Yamaguchi2020,He2024} is the most straightforward and computationally efficient semi-implicit approach among the others. Recent research, based on the incremental fractional step method \citep{Kularathna2021, Yuan2023}, improved its performance in geotechnical problems by removing the permeability restriction to the critical time step in the traditional fractional step method. A stabilised semi-implicit double-point incremental fractional step MPM formulation has been recently developed by the Authors \citep{Xie2025}, based on the conventional single-point approach proposed by \citet{Kularathna2021, Yuan2023}. It was concluded that, employing two distinct point sets for soil and water, markedly improves accuracy and robustness, including when benchmarked against single-point counterparts, using the same time-integration framework. However, applying that method to \emph{intense} soil–water interaction, such as dam-break impacts, rapid seepage across sharp porosity transitions, and landslide-induced tsunamis, still requires substantial developments to ensure stability and accuracy.

{Several double-point MPMs have recently been developed for soil--water interaction, using either semi-implicit fractional step \citep{Yamaguchi2020, He2024} or monolithic \citep{Chandra2024b} approaches, each with distinct strengths. However, applying these formulations to the full spectrum of problems addressed here, including pure free-surface flow with wave breaking, rapid seepage across sharp porosity transitions, and landslide-generated tsunamis, poses specific challenges that motivate the present developments. The traditional fractional step method adopted by \citet{Yamaguchi2020} and \citet{He2024} imposes a permeability-dependent time-step restriction; for practical geotechnical problems with relatively low permeability, this necessitates very small time steps, increasing computational cost considerably. In addition, as reported in \citet{Chandra2024b}, the traditional fractional step method suffers from a volume conservation issue, whereby the water volume gradually decreases over the course of the simulation. The monolithic approach of \citet{Chandra2024b} avoids these problems but requires solving a fully coupled system at each step, also increasing the computational cost. The MPM-FEM hybrid of \citet{Pan2021} treats the fluid phase within an Eulerian FEM framework on the background mesh rather than tracking it with Lagrangian material points, which may introduce numerical diffusion at sharp free-surface interfaces. Furthermore, existing double-point MPM studies of landslide-generated tsunamis and intense soil-water interaction have typically relied on the Drucker-Prager constitutive model \citep{He2024, Pan2021}, which does not capture volume changes or density evolution in granular materials. The present work aims to address these limitations by extending the incremental fractional step method, which removes both the permeability-dependent time-step restriction and the volume conservation problem, from confined pore-water problems to free-surface flows, whilst incorporating an advanced critical-state constitutive model.}

The objective of this paper is to develop a new stabilised semi-implicit, two-phase, double-point Material Point Method (MPM) \citep{Xie2023b} to reliably model problems that involve free water and strong two-way interaction with deformable porous beds. The contributions are fivefold: (i) {for the first time, we incorporate} a Forchheimer-type non-linear extension of Darcy's law within an incremental fractional-step algorithm without enlarging the implicit solve, enabling high–Reynolds-number interphase momentum exchange; (ii) we adopt a hyperelastic Newtonian constitutive model for slightly compressible, viscous water in a semi-implicit coupled MPM, improving numerical stability in fast free-surface motion and wave breaking; (iii) we develop a mixed Taylor Particle-in-Cell scheme for B-spline MPM that curbs spurious velocity oscillations with minimal dissipation;  (iv) we introduce interface-robust ingredients tailored to higher-order spline bases, including a nodal free-surface detection for pressure boundary enforcement and a smoothly varying porosity–permeability map based on the Kozeny–Carman relation; and (v) to address quadrature errors when water points traverse sharp porosity jumps, we employ a one-step $\delta$-correction to maintain nodal volume consistency.

The remainder of this paper begins with a review of the existing semi-implicit formulation in Sect.~\ref{sec:baseline} to provide a baseline for further developments towards modelling soil-water coupled problems with free-surface flows, which are presented in Sect.~\ref{sec:improvements}. We validate the enhanced formulation against experimental benchmarks spanning pure-water dam break \citep{Lobovsky2014}, dam-break seepage through a porous barrier \citep{Liu1999}, two granular-collapse tsunami experiments \citep{Rauter2022,Sarlin2022}, {and a dam-break wave over a movable granular bed \citep{Spinewine2005}} in Sect.~\ref{sec:num_examples}. For the granular phase we use an advanced NorSand model in a large-strain framework \citep{Borja2006,Jefferies2015}, calibrated against triaxial data for glass beads \citep{Wu2017}, to highlight the importance of density evolution and dilatancy compared with the simpler Drucker–Prager or $\mu$–$I$ rheological laws used in previous studies \citep{Pan2021,Rauter2022,He2024}. The results show that the proposed semi-implicit double-point MPM accurately reproduces free-surface evolution, pressure time histories, seepage fronts, wave gauges{, and erosion} for a wide range of soil-water interaction problems, while remaining stable even when interactions are intense.

\section{Governing equations and semi-implicit coupled formulation}\label{sec:baseline}
This section summarises the semi-implicit coupled formulation based on the incremental fractional step method to introduce symbols and notation and to provide
a starting point for improving this method towards modelling free-surface flows.

\subsection{Conservation laws for the mixture}

A fully saturated soil–water mixture with separate Lagrangian descriptions for soil (\(s\)) and water (\(w\)) has been considered.
The material points for each phase map refer to spatial coordinates via
\(\boldsymbol{x}_\alpha=\boldsymbol{\chi}_\alpha(\boldsymbol{X}_\alpha,t)\),
with velocities \(\boldsymbol{v}_\alpha=\mathrm{d}\boldsymbol{x}_\alpha/\mathrm{d}t\) and accelerations \(\boldsymbol{a}_\alpha=\mathrm{d}\boldsymbol{v}_\alpha/\mathrm{d}t\), \(\alpha\in\{s,w\}\).
Both the soil's grains and water are treated as incompressible; porosity \(n\) evolves in time.

Following the classical mixture-theory framework for porous media
\citep{Morland1992,deBoerEhlers1990} and its recent use in debris-flow
modelling by \citet{Meng2022}, the water volume fraction is denoted by \(n\),
while the solid volume fraction is \(1-n\).
The corresponding partial densities are
\begin{equation}
	\rho_s^{p}=(1-n)\rho_s,
	\qquad
	\rho_w^{p}=n\rho_w .
	\label{eq:partial-density}
\end{equation}

For incompressible soil grains and incompressible water, the phase mass balance
equations are
\begin{subequations}\label{eq:mass-mixture}
	\begin{align}
		\frac{\text{d}_s(1-n)}{\text{d}t}+(1-n)\nabla\cdot\boldsymbol{v}_s &= 0,\\
		\frac{\text{d}_w n}{\text{d}t}+n\nabla\cdot\boldsymbol{v}_w &= 0,
	\end{align}
\end{subequations}
where \(\text{d}_s/\text{d}t\) and \(\text{d}_w/\text{d}t\) denote material derivatives following the soil and
water phases, respectively. Combining the two phase mass balance equations gives
\begin{equation}
	(1-n)\nabla\cdot\boldsymbol{v}_s+n\nabla\cdot\boldsymbol{v}_w=0 .
	\label{eq:mixture-div}
\end{equation}

For phase \(\alpha\in\{s,w\}\), the local momentum balance can be written in
conservative form as
\begin{equation}
	\frac{\partial}{\partial t}
	\left(\rho_\alpha^{p}\boldsymbol{v}_\alpha\right)
	+
	\nabla\cdot
	\left(
		\rho_\alpha^{p}
		\boldsymbol{v}_\alpha\otimes\boldsymbol{v}_\alpha
	\right)
	=
	\nabla\cdot\boldsymbol{\sigma}_\alpha
	+
	\rho_\alpha^{p}\boldsymbol{b}
	+
	\boldsymbol{\beta}_\alpha ,
	\label{eq:phase-mom-cons}
\end{equation}
where \(\boldsymbol{\sigma}_\alpha\) is the partial stress of phase
\(\alpha\), and \(\boldsymbol{\beta}_\alpha\) is the interaction force exerted
on phase \(\alpha\) by the other phase. The interaction forces satisfy
\begin{equation}
	\boldsymbol{\beta}_s+\boldsymbol{\beta}_w=\boldsymbol{0}.
	\label{eq:interaction-balance}
\end{equation}
Using the phase mass balance in Eq.~\eqref{eq:mass-mixture}, Eq.~\eqref{eq:phase-mom-cons} is equivalently written in
material-acceleration form as
\begin{equation}
	\rho_\alpha^{p}\boldsymbol{a}_\alpha
	=
	\nabla\cdot\boldsymbol{\sigma}_\alpha
	+
	\rho_\alpha^{p}\boldsymbol{b}
	+
	\boldsymbol{\beta}_\alpha .
	\label{eq:phase-mom-material}
\end{equation}

Following Terzaghi's effective-stress principle \citep{Terzaghi1943} and the
partial-stress decomposition commonly used in mixture-theory formulations of
porous media \citep{deBoerEhlers1990,Meng2022}, the partial stresses are
decomposed as
\begin{subequations}\label{eq:partial-stresses}
	\begin{align}
		\boldsymbol{\sigma}_s
		&=
		\boldsymbol{\sigma}'
		+
		(1-n)p_w\boldsymbol{I},
		\label{eq:partial-stress-s}
		\\
		\boldsymbol{\sigma}_w
		&=
		n p_w\boldsymbol{I}
		+
		\boldsymbol{\tau}_w ,
		\label{eq:partial-stress-w}
	\end{align}
\end{subequations}
where \(\boldsymbol{\sigma}'\) is the soil effective stress tensor, \(p_w\) is the water pressure and \(\boldsymbol{I}\) is the second-order identity tensor. Here \(\boldsymbol{\tau}_w\) is the partial non-pressure (deviatoric) stress carried by the
water phase. In the hydrostatic baseline formulation,
\(\boldsymbol{\tau}_w=\boldsymbol{0}\). For a Newtonian viscous water phase,
\(\boldsymbol{\tau}_w\) may be written as
\begin{equation}
	\boldsymbol{\tau}_w
	=
	n
	\left[
		2\mu\boldsymbol{D}_w
		-
		\frac{2}{3}\mu\,\mathrm{tr}(\boldsymbol{D}_w)\boldsymbol{I}
	\right],
	\qquad
	\boldsymbol{D}_w
	=
	\frac{1}{2}
	\left[
		\nabla\boldsymbol{v}_w+
		(\nabla\boldsymbol{v}_w)^{T}
	\right],
	\label{eq:tauw-newtonian}
\end{equation}
where $\mu = 0.001$ Pa$\cdot$s is the dynamic viscosity of water at 20 °C. When the hyperelastic water model is used in
Sect.~\ref{sec:hyperelastic_water}, the same term is evaluated from the
non-pressure part of the intrinsic water Cauchy stress,
\begin{equation}
	\boldsymbol{\tau}_w
	=
	n\left(
		\boldsymbol{\sigma}^{\mathrm{HE}}_w
		-
		p_w\boldsymbol{I}
	\right).
\end{equation}

The pressure parts of Eq.~\eqref{eq:partial-stresses} give
\begin{subequations}\label{eq:pressure-product-rule}
	\begin{align}
		\nabla\cdot\left[(1-n)p_w\boldsymbol{I}\right]
		&=
		(1-n)\nabla p_w
		-
		p_w\nabla n,
		\label{eq:pressure-product-s}
		\\
		\nabla\cdot\left[n p_w\boldsymbol{I}\right]
		&=
		n\nabla p_w
		+
		p_w\nabla n .
		\label{eq:pressure-product-w}
	\end{align}
\end{subequations}
Following the mixture-theory interpretation that Darcy's law provides a
constitutive closure for the interaction body force in porous media
\citep{Morland1992}, and following the pressure-interaction treatment of
\citet{Meng2022}, the porosity-gradient pressure contribution is included in
the interphase pressure interaction. The interaction force is
therefore decomposed as
\begin{subequations}\label{eq:interaction-force}
	\begin{align}
		\boldsymbol{\beta}_s
		&=
		p_w\nabla n
		+
		\hat{\boldsymbol{p}},
		\label{eq:interaction-force-s}
		\\
		\boldsymbol{\beta}_w
		&=
		-
		p_w\nabla n
		-
		\hat{\boldsymbol{p}},
		\label{eq:interaction-force-w}
	\end{align}
\end{subequations}
where \(\hat{\boldsymbol{p}}\) denotes the Darcy or Forchheimer drag force
associated with the relative motion between the soil and water phases.

Substituting Eqs.~\eqref{eq:partial-stress-s} and
\eqref{eq:interaction-force-s} into Eq.~\eqref{eq:phase-mom-material}, and using
Eq.~\eqref{eq:pressure-product-s}, gives the soil phase momentum equation as
\begin{align}
	\rho_s^{p}\boldsymbol{a}_s
	&=
	\nabla\cdot\boldsymbol{\sigma}'
	+
	(1-n)\nabla p_w
	-
	p_w\nabla n
	+
	\rho_s^{p}\boldsymbol{b}
	+
	p_w\nabla n
	+
	\hat{\boldsymbol{p}}
	\notag\\
	&=
	\nabla\cdot\boldsymbol{\sigma}'
	+
	(1-n)\nabla p_w
	+
	\rho_s^{p}\boldsymbol{b}
	+
	\hat{\boldsymbol{p}}.
	\label{eq:mom-s}
\end{align}
Similarly, substituting Eqs.~\eqref{eq:partial-stress-w} and
\eqref{eq:interaction-force-w} into Eq.~\eqref{eq:phase-mom-material}, and
using Eq.~\eqref{eq:pressure-product-w}, gives
\begin{align}
	\rho_w^{p}\boldsymbol{a}_w
	&=
	n\nabla p_w
	+
	p_w\nabla n
	+
	\nabla\cdot\boldsymbol{\tau}_w
	+
	\rho_w^{p}\boldsymbol{b}
	-
	p_w\nabla n
	-
	\hat{\boldsymbol{p}}
	\notag\\
	&=
	n\nabla p_w
	+
	\nabla\cdot\boldsymbol{\tau}_w
	+
	\rho_w^{p}\boldsymbol{b}
	-
	\hat{\boldsymbol{p}}.
	\label{eq:mom-w}
\end{align}

Therefore, the final phase momentum equations are
\begin{subequations}\label{eq:momentum}
	\begin{align}
		\rho_s^{p}\boldsymbol{a}_s
		&=
		\nabla\cdot\boldsymbol{\sigma}'
		+
		(1-n)\nabla p_w
		+
		\rho_s^{p}\boldsymbol{b}
		+
		\hat{\boldsymbol{p}},
		\label{eq:momentum-s}
		\\
		\rho_w^{p}\boldsymbol{a}_w
		&=
		n\nabla p_w
		+
		\nabla\cdot\boldsymbol{\tau}_w
		+
		\rho_w^{p}\boldsymbol{b}
		-
		\hat{\boldsymbol{p}}.
		\label{eq:momentum-w}
	\end{align}
\end{subequations}
The term \((1-n)\nabla p_w\) represents the buoyancy
contribution associated with the intrinsic pore-water pressure gradient acting
on the soil skeleton. The porosity-gradient pressure term \(p_w\nabla n\) is
not neglected; instead, it is included in the interphase pressure interaction
and cancels from the final phase momentum equations. The non-pressure water
stress contributes through \(\nabla\cdot\boldsymbol{\tau}_w\), consistent with
the water-phase momentum equation of \citet{Meng2022}.

For the linear Darcy interaction, consistent with Darcy-type flow through porous media \citep{Bear1972,Morland1992}, the drag term is written as
\begin{equation}
	\hat{\boldsymbol{p}}
	=
	\frac{n^2\rho_w g}{k}
	\left(
		\boldsymbol{v}_w-\boldsymbol{v}_s
	\right),
	\label{eq:darcy}
\end{equation}
where \(k\) is the hydraulic conductivity and \(g\) is the gravitational
acceleration.

Adding Eqs.~\eqref{eq:momentum-s} and \eqref{eq:momentum-w} gives the mixture
momentum equation
\begin{equation}
	\rho_s^{p}\boldsymbol{a}_s
	+
	\rho_w^{p}\boldsymbol{a}_w
	=
	\nabla\cdot\boldsymbol{\sigma}'
	+
	\nabla p_w
	+
	\nabla\cdot\boldsymbol{\tau}_w
	+
	\rho_s^{p}\boldsymbol{b}
	+
	\rho_w^{p}\boldsymbol{b}.
	\label{eq:mixture-mom}
\end{equation}
Equation~\eqref{eq:mixture-mom} is consistent with the mixture total stress
\begin{equation}
	\boldsymbol{\sigma}
	=
	\boldsymbol{\sigma}'
	+
	p_w\boldsymbol{I}
	+
	\boldsymbol{\tau}_w .
	\label{eq:mixture-total-stress}
\end{equation}

\subsection{Semi-implicit incremental fractional-step scheme}
The idea of a semi-implicit incremental fractional-step scheme is to introduce
intermediate predictor velocities \(\boldsymbol{v}_\alpha^*\) and split the acceleration as
\begin{equation}
	\boldsymbol{a}_\alpha^{\,t+1}
	=
	\frac{\boldsymbol{v}_\alpha^{\,t+1}-\boldsymbol{v}_\alpha^t}{\Delta t}
	=
	\frac{\boldsymbol{v}_\alpha^{\,t+1}-\boldsymbol{v}_\alpha^*}{\Delta t}
	+
	\frac{\boldsymbol{v}_\alpha^*-\boldsymbol{v}_\alpha^t}{\Delta t}.
	\label{eq:accel-split}
\end{equation}
The pressure increment is solved implicitly, while the effective soil stress, the known
pressure contribution, the non-pressure water stress, body force and interphase drag are
included in the predictor step. In particular, \(\boldsymbol{\tau}_w\) is evaluated explicitly
from known fields at time \(t\), or from the previously updated water constitutive state,
so that no additional unknown is introduced into the pressure projection.

\paragraph{Predictor (explicit) step:}
Using known fields at \(t\), update the intermediate water velocity
\(\boldsymbol{v}_w^*\) from
\begin{subequations}\label{eq:water-predictor}
	\begin{align}
		\rho_w^{p}
		\left(
			\frac{\boldsymbol{v}_w^*-\boldsymbol{v}_w^{\,t}}{\Delta t}
		\right)
		&=
		n\nabla p_w^{\,t}
		+
		{ \nabla\cdot\boldsymbol{\tau}_w^{\,t} }
		+
		\rho_w^{p}\boldsymbol{b}^{\,t+1}
		-
		\hat{\boldsymbol{p}},
		\label{eq:water-pred-a}
		\\
		\hat{\boldsymbol{p}}
		&=
		\frac{n^2\rho_w g}{k}
		\left(
			\boldsymbol{v}_w^*-\boldsymbol{v}_s^{\,t}
		\right).
		\label{eq:water-pred-b}
	\end{align}
\end{subequations}
The term \(p_w^{\,t}\nabla n\) does not appear in
Eq.~\eqref{eq:water-pred-a}, because it is included in the interphase pressure
interaction and cancels from the final phase momentum equation.

The intermediate soil velocity \(\boldsymbol{v}_s^*\) is then obtained from the mixture
predictor:
\begin{align}
	\rho_s^{p}
	\left(
		\frac{\boldsymbol{v}_s^*-\boldsymbol{v}_s^{\,t}}{\Delta t}
	\right)
	+
	\rho_w^{p}
	\left(
		\frac{\boldsymbol{v}_w^*-\boldsymbol{v}_w^{\,t}}{\Delta t}
	\right)
	&=
	\nabla\cdot\boldsymbol{\sigma}'^{\,t}
	+
	\nabla p_w^{\,t}
	+
	{ \nabla\cdot\boldsymbol{\tau}_w^{\,t} }
	+
	\rho_s^{p}\boldsymbol{b}^{\,t+1}
	+
	\rho_w^{p}\boldsymbol{b}^{\,t+1}.
	\label{eq:mixture-predictor}
\end{align}

\paragraph{Pressure projection (implicit) step:}
The pressure increment
\(\delta p_w^{\,t+1}=p_w^{\,t+1}-p_w^{\,t}\) is solved by enforcing the mixture
divergence constraint at \(t+1\). Since the non-pressure water stress has already been
included explicitly in the predictor step, the projection equation remains
\begin{equation}
	(1-n)\nabla\cdot\boldsymbol{v}_s^*
	+
	n\nabla\cdot\boldsymbol{v}_w^*
	+
	\Delta t
	\left(
		\frac{1-n}{\rho_s}
		+
		\frac{n}{\rho_w}
	\right)
	\nabla^2\delta p_w^{\,t+1}
	=
	0 .
	\label{eq:pressure-proj}
\end{equation}
Zero pressure is imposed on the detected free-water surface.

\paragraph{Velocity corrector step:}
The phase velocities are corrected using only the pressure increment:
\begin{subequations}\label{eq:vel-correct}
	\begin{align}
		\rho_w^{p}
		\left(
			\frac{\boldsymbol{v}_w^{\,t+1}-\boldsymbol{v}_w^*}{\Delta t}
		\right)
		&=
		n\nabla\delta p_w^{\,t+1},
		\label{eq:vel-correct-w}
		\\
		\rho_s^{p}
		\left(
			\frac{\boldsymbol{v}_s^{\,t+1}-\boldsymbol{v}_s^*}{\Delta t}
		\right)
		+
		\rho_w^{p}
		\left(
			\frac{\boldsymbol{v}_w^{\,t+1}-\boldsymbol{v}_w^*}{\Delta t}
		\right)
		&=
		\nabla\delta p_w^{\,t+1}.
		\label{eq:vel-correct-s}
	\end{align}
\end{subequations}
After the velocity correction, the water constitutive state is updated {and
\(\boldsymbol{\tau}_w^{\,t+1}\) is evaluated} for use in the next predictor step.

The weak forms of Eqs.~\eqref{eq:water-predictor}–\eqref{eq:vel-correct} are evaluated over the current soil and water domains using a standard Galerkin approach with material points as quadrature points. Soil and water point data are mapped to common grid nodes (equal-order interpolation for displacement/velocity and pressure). Detailed discrete expressions follow \citet{Xie2025} and are omitted here for brevity.

\section{Towards modelling free-surface flows}\label{sec:improvements}
\subsection{Non-linear Darcy law}\label{sec:non_linear_darcy}

In the soil-water coupled formulations \citep{Bandara2015,Kularathna2021,Yuan2023,Xie2025}, the flow in the pore is usually assumed to follow the linear Darcy's law. In problems involving free water, such as landslide-induced tsunamis, this assumption may no longer be valid because the Reynolds number can reach very high values. To address this issue, the momentum interaction ($\hat{\pmb{p}}$) between the soil and water phase under the linear Darcy law shown in Eq.~\eqref{eq:darcy} is usually replaced by the Forchheimer-like equation \citep{Chandra2024b,He2024,Pan2021}
\begin{equation}\label{eq:nonlinearDarcy}
	\hat{\pmb{p}} = \frac{n^2\,\rho_w\,g}{k}\cdot(\pmb{v}_w - \pmb{v}_s) + \frac{1.75n\,{\rho}_w\sqrt{n\,\rho_w\,g}\cdot\abs{\pmb{v}_w - \pmb{v}_s}}{\sqrt{150k\,\mu}} \cdot(\pmb{v}_w - \pmb{v}_s) \,.
\end{equation}
Compared with Eq.~\eqref{eq:darcy}, a non-linear term related to the absolute velocity difference between the soil and water phases is added in Eq.~\eqref{eq:nonlinearDarcy}. In this research, the viscosity of water is assumed to be a constant value for both pore and free water.

To our best knowledge, it is the first time that the non-linear Darcy law has been implemented in an incremental fractional step method. To implement Eq.~\eqref{eq:nonlinearDarcy} into the incremental fractional step method without additional computational cost, we use the velocity of water and soil from the last time step ($\pmb{v}^t_w$ and $\pmb{v}^t_s$) to calculate the absolute value of the second term of Eq.~(\ref{eq:nonlinearDarcy}): 
\begin{equation}\label{eq:momentum_term_new}
	\hat{\pmb{p}} = \frac{n^2\,\rho_w\,g}{k}\cdot(\pmb{v}^{*}_w - \pmb{v}^t_s) + \frac{1.75n\,{\rho}_w\sqrt{n\,\rho_w\,g}\cdot\abs{\pmb{v}^t_w - \pmb{v}^t_s}}{\sqrt{150k\,\mu}}\cdot(\pmb{v}^{*}_w - \pmb{v}^t_s)  \,.
\end{equation}
In this case, the intermediate water velocity $\pmb{v}^{*}_w$ is still the only unknown value in Eqs.~\eqref{eq:momentum_term_new} and \eqref{eq:water-pred-a}, retaining the advantage of efficient numerical implementation for the proposed formulation. {Alternatively, Eq.~\eqref{eq:nonlinearDarcy} could be implemented by using both intermediate velocities $\boldsymbol{v}_w^{*}$ and $\boldsymbol{v}_s^{*}$ in all terms. In this case, however, Eqs.~\eqref{eq:water-predictor}-\eqref{eq:pressure-proj} would need to be solved simultaneously, coupling the predictor and pressure projection steps and significantly increasing the computational cost. The present linearisation, using $\boldsymbol{v}_w^t$ and $\boldsymbol{v}_s^t$ from the previous time step, preserves the sequential structure of the incremental fractional step algorithm: the intermediate water velocity $\boldsymbol{v}_w^{*}$ remains the only unknown in Eq.~\eqref{eq:momentum_term_new}, and the overall efficiency of the method is maintained. The numerical examples in Sect.~\ref{sec:num_examples} confirm that this linearisation yields stable and accurate results for high-Reynolds-number free-surface flows, demonstrating the robustness of this approach.}

\subsection{Hyperelastic constitutive model for water}\label{sec:hyperelastic_water}

In the baseline formulation introduced in Sect.~\ref{sec:baseline}, the water contribution is represented by the classical hydrostatic pressure term \(p_w\boldsymbol{I}\), corresponding to an incompressible pore-water pressure field. This assumption is commonly adopted in existing semi-implicit coupled MPMs \citep{Yamaguchi2020,Kularathna2021,Pan2021,Yuan2023,He2024,Chandra2024b,Xie2025}. In reality, water is slightly compressible and possesses finite viscosity. These features cannot be accurately represented by the simple hydrostatic term \(p_w\boldsymbol{I}\). Although this approximation is reasonable for many pore-water problems, it may not be sufficiently robust for fast-moving free-water flows. {Therefore, in the enhanced formulation, the pressure part \(p_w\boldsymbol{I}\) is retained in the semi-implicit pressure projection, while the compressibility and viscous contributions of the water constitutive model are introduced through the non-pressure water stress \(\boldsymbol{\tau}_w\) defined in Eq.~\eqref{eq:tauw-hyperelastic}.}

\citet{Molinos2023} applied this model in MPM to study fluid-structure interaction problems. {The intrinsic Cauchy stress tensor for slightly compressible water is denoted by \(\boldsymbol{\sigma}^{\mathrm{HE}}_w\) and is obtained as \citep{Molinos2023}:}
\begin{equation}\label{eq:sig_w}
	\boldsymbol{\sigma}^{\mathrm{HE}}_w
	=
	p_w\boldsymbol{I}
	-
	\frac{\kappa_f}{m}
	\left(
		J^{-m}-1
	\right)
	\boldsymbol{I}
	-
	\frac{2}{3}\mu\,\mathrm{tr}(\boldsymbol{D})\boldsymbol{I}
	+
	2\mu\boldsymbol{D}.
\end{equation}
where
\begin{equation}
	\boldsymbol{D}
	=
	\mathrm{sym}
	\left(
		\dot{\boldsymbol{F}}\boldsymbol{F}^{-1}
	\right)
	\label{eq:water-D}
\end{equation}
is the spatial rate of deformation tensor; \(\boldsymbol{F}\) is the deformation gradient tensor; \(\dot{\boldsymbol{F}}=\nabla\boldsymbol{v}_w\) is the material velocity gradient of the water phase; \(J=\det(\boldsymbol{F})\) is the Jacobian of the deformation gradient; \(m\) is a polytropic index, taken as \(7\) to reproduce slightly compressible fluids; and $\kappa_f$ is the compressibility of fluid which equals $2.2 \times 10^6$ Pa for water \citep{Molinos2023}. 

{Consistent with the phase-stress decomposition in Eq.~\eqref{eq:partial-stress-w}, the water partial stress is written as}
\begin{equation}
	{
	\boldsymbol{\sigma}_w
	=
	n p_w\boldsymbol{I}
	+
	\boldsymbol{\tau}_w ,
	}
	\label{eq:water-partial-stress-he}
\end{equation}
{where the partial non-pressure water stress is evaluated from the non-pressure part of Eq.~\eqref{eq:sig_w}:}
\begin{equation}
	{
	\boldsymbol{\tau}_w
	=
	n
	\left(
		\boldsymbol{\sigma}^{\mathrm{HE}}_w
		-
		p_w\boldsymbol{I}
	\right)
	=
	n
	\left[
		-
		\frac{\kappa_f}{m}
		\left(
			J^{-m}-1
		\right)
		\boldsymbol{I}
		-
		\frac{2}{3}\mu\,\mathrm{tr}(\boldsymbol{D})\boldsymbol{I}
		+
		2\mu\boldsymbol{D}
	\right].
	}
	\label{eq:tauw-hyperelastic}
\end{equation}
{Substituting Eq.~\eqref{eq:tauw-hyperelastic} into the water momentum equation Eq.~\eqref{eq:momentum-w} gives the additional term \(\nabla\cdot\boldsymbol{\tau}_w\). In the semi-implicit fractional-step scheme, this term is evaluated explicitly in the predictor step, as shown in Eqs.~\eqref{eq:water-pred-a} and \eqref{eq:mixture-predictor}. Therefore, the implicit pressure projection in Eq.~\eqref{eq:pressure-proj} remains a scalar pressure-increment solve and the size of the implicit system is not increased.}

{The baseline hydrostatic formulation is recovered by neglecting the compressibility and viscous shear contributions in Eq.~\eqref{eq:sig_w}, which gives \(\boldsymbol{\sigma}^{\mathrm{HE}}_w=p_w\boldsymbol{I}\) and therefore \(\boldsymbol{\tau}_w=\boldsymbol{0}\) in Eq.~\eqref{eq:tauw-hyperelastic}.}

This is the first implementation of a water constitutive model in a semi-implicit coupled MPM, to the best of our knowledge. The motivation for adopting a hyperelastic, slightly compressible model is twofold. First, water at room temperature is often treated as perfectly incompressible, but it is in fact slightly compressible and possesses a finite dynamic viscosity. The classical hydrostatic representation \(p_w\boldsymbol{I}\) used in existing semi-implicit coupled MPMs neglects both features. {Adopting Eqs.~\eqref{eq:sig_w} and \eqref{eq:tauw-hyperelastic} provides a more faithful description of the non-pressure water response while retaining the pressure field \(p_w\) as the primary unknown in the semi-implicit projection.} Second, from a numerical perspective, the slight compressibility regularises the strict incompressibility constraint, which in semi-implicit MPM can cause spurious pressure oscillations near free surfaces and material point clumping during fast flows. The high value of \(\kappa_f\) keeps volume changes very small, hence the macroscopic behaviour remains close to that of incompressible water, while it provides physical damping of high-frequency pressure modes and contributes to stable solutions. The viscous shear term provides physical dissipation that replaces artificial numerical damping for fast free-surface motion. This combination is particularly beneficial for high-Reynolds-number flows such as wave breaking and landslide impact, as demonstrated in Sect.~\ref{sec:dambreak}, where simulations without the hyperelastic water model diverge under intense free-surface motion, while the proposed formulation remains stable at negligible additional computational cost.

\subsection{A stable motion integration strategy}\label{sec:1pTPIC}
The FLuid Implicit Particle (FLIP) method developed by \citet{Brackbill1986} is frequently used in MPM to update the material points' velocity. The FLIP method was developed by \citet{Brackbill1986} to improve the numerical dissipation in the Particle-in-Cell (PIC) scheme proposed by \citet{Harlow1963}, as the latter method had limited applications due to its massive energy dissipation. In MPM under the FLIP scheme, the material point velocity is updated by
\begin{equation}\label{eq:FLIP}
    \pmb{v}_{p}^{t + 1} = \pmb{v}_{p}^{t} + \sum_{I=1}^{{N}_{n}} {S}_{Ip} \, \Delta\pmb{v}_{I}^{t + 1} \,,
\end{equation}
where ${N}_{n}$ is the number of grid nodes that influence the material point, and ${S}_{Ip}$ is the shape function associated with node $I$ evaluated at the position of material point $p$. Under the PIC scheme, the material point velocity is updated by mapping the grid velocity directly to the material point:
\begin{equation}\label{eq:PIC}
    \pmb{v}_{p}^{t + 1} = \sum_{I=1}^{{N}_{n}} {S}_{Ip} \, \pmb{v}_{I}^{t + 1} \,.
\end{equation}

The FLIP method has been shown to have no energy dissipation in MPM \citep{Duverger2024}. However, it results in a significant oscillation in the velocity field. \citet{Nakamura2023} reported examples of this oscillation in the FLIP scheme. Whilst using the FLIP scheme to update the soil or pore water's velocity typically results in a stable displacement field despite velocity oscillation, significant instabilities in the displacement field can be observed in MPM with the FLIP method for free-water problems involving high Reynolds numbers, such as the dam break problem. This can result in portions of the water material points dispersing throughout the computational domain without returning to the main flow \citep{Chandra2024b}.

To overcome this instability, the Affine PIC (APIC) \citep{Jiang2015} and Taylor PIC (TPIC) \citep{Nakamura2023} were proposed, under the PIC family. \citet{Nakamura2023} developed TPIC following the Taylor FLIP (TFLIP) method proposed by \citet{Wallstedt2007}. According to \citet{Nakamura2023} and \citet{Chandra2024b}, APIC and TPIC demonstrate very similar performance, while the latter method is mathematically simpler and easier to implement. In APIC, TPIC and TFLIP, at the beginning of each computational cycle, the grid velocity $\pmb{v}_{I}^t$ is obtained from an affine velocity $\Tilde{\pmb{v}}^t_p$:
\begin{equation}
\pmb{v}_{I}^t = \frac{1}{\pmb{M}_{I} } \sum_{p=1}^{{N}_{p}} {S}_{Ip} \, m_{p} \, \Tilde{\pmb{v}}^t_p(\pmb{x}_I) \,.
\end{equation}
The affine velocity $\Tilde{\pmb{v}}^t_p(\pmb{x}_I)$ at grid $\pmb{x}_I$ can be obtained by taking the first-order Taylor series approximation about the location of the material point $\pmb{x}^t_p$:
\begin{equation}\label{eq:TPIC_v}
    \Tilde{\pmb{v}}^t_p(\pmb{x}_I) = \pmb{v}^t_p + \nabla\pmb{v}^t_p \cdot (\pmb{x}_I - \pmb{x}^t_p) \,,
\end{equation}
resulting in TFLIP or TPIC depending on the mapping schemes (Eq.~(\ref{eq:FLIP}) or Eq.~(\ref{eq:PIC})) used to update the material point velocity. The APIC follows a different approach to obtaining the affine velocity (see \citet{Jiang2015} and \citet{Nakamura2023} for details).

Although APIC and TPIC show promising performance in mitigating velocity oscillation, their dissipations are as high as in the original PIC \citep{Duverger2024}. Therefore, extra care must be taken when applying a method from the PIC family, including APIC and TPIC, to an engineering problem. Because of that, to stabilise the velocity field of the water phase whilst conserving most of the energy, we update the water material point velocity using a linear combination of the TFLIP and TPIC methods:
\begin{equation}\label{eq:mFLIP}
    \pmb{v}_{p}^{t + 1} = \beta\sum_{I=1}^{{N}_{n}} {S}_{Ip}\, \pmb{v}_{I}^{t + 1} + (1-\beta)(\pmb{v}_{p}^{t} + \sum_{I=1}^{{N}_{n}} {S}_{Ip} \, \Delta\pmb{v}_{I}^{t + 1}) \,,
\end{equation}
where $\beta$ is a number between 0 and 1 that controls the contributions of the TFLIP scheme. In this research, only 1\% of TPIC (that is, $\beta = 0.01$) is used to maintain most of the energy, denoted as the TPIC (1\%) method. 

In TPIC-based MPM research \citep{Nakamura2023,Chandra2024b}, a B-spline shape function is typically used together with a weighted least squares (WLS) kernel correction. Without this, the algorithm can easily diverge because the interpolation error near the boundaries in a B-spline MPM is magnified by the affine velocity $\Tilde{\pmb{v}}^t_p(\pmb{x}_I)$ \citep{Nakamura2023,Chandra2024b}. However, this boundary error occurs due to the use of ghost nodes (i.e., nodes outside the computational domain) in their implementation of an explicitly defined B-spline shape function. This error can be completely avoided by using a B-spline formed by open knot vectors using the Cox-de Boor formula as mentioned in \citet{Xie2023b}. Therefore, in this case, a kernel correction becomes unnecessary. Furthermore, an explicitly defined B-spline shape function cannot ensure a positive definite mass matrix. In \citet{Chandra2024b}, negative mass is avoided by applying kernel correction iteratively, resulting in additional computational costs. In contrast, a B-spline shape function formed by open knot vectors is positive everywhere, ensuring a positive definite mass matrix.

\begin{figure}[htbp]
	\begin{center}
		\includegraphics[width=0.7\textwidth]{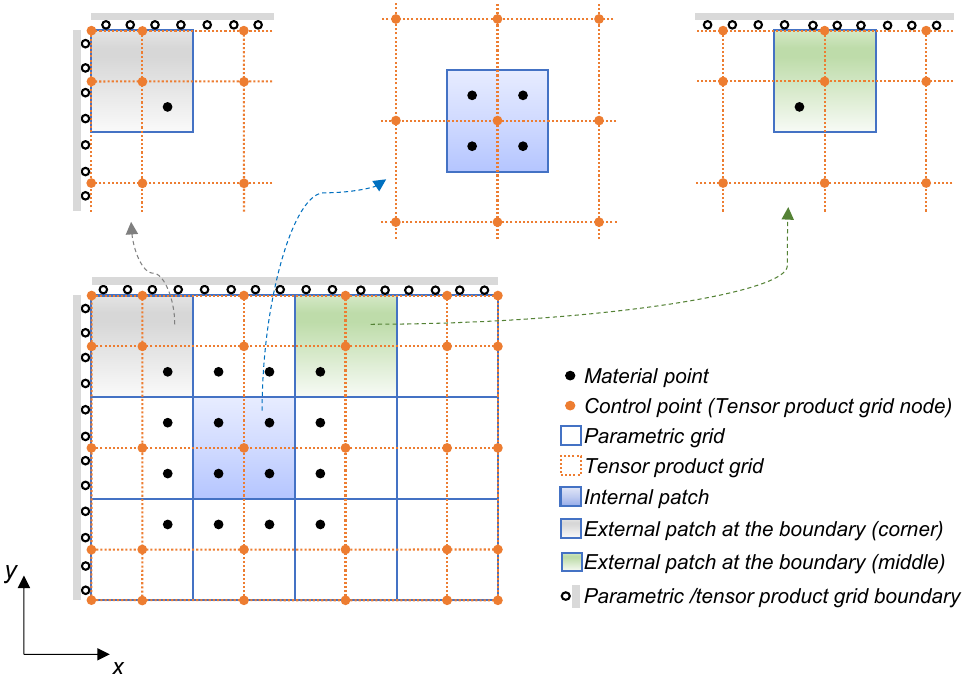}
		\caption{Spatial relationships between material point, parametric grid and 2-D quadratic B-spline grid formed by the tensor product of open knot vectors}
		\label{QuadraticBsplineGrids}
	\end{center}
\end{figure}

As shown in Eq.~(\ref{eq:TPIC_v}), to calculate the affine velocity $\Tilde{\pmb{v}}^t_p(\pmb{x}_I)$, the spatial relationship between the grid nodes and the material points is needed. In a B-spline MPM formed by open knot vectors, calculations are carried out in the tensor product grid nodes, which have different coordinates than the parametric grid nodes \citep{Xie2023b}. However, the spatial relationship between the material points and the nodes of the tensor product grid has not yet been illustrated. Figure.~\ref{QuadraticBsplineGrids} presents the spatial relationship for a quadratic B-spline MPM \citep{Cottrell2009}. In a quadratic B-spline MPM, we have three possible spatial arrangements of the tensor product grid nodes depending on the location of their controlled material point, as shown in Fig.~\ref{QuadraticBsplineGrids}. Therefore, the coordinate differences between the nodes of the tensor product grid and the associated material points should be calculated accordingly.

\subsection{Smooth soil-water interface}\label{sec:sw_interface}

In some soil-water coupled MPMs \citep{Bandara2015, Liu2017, Xie2025}, a sharp soil-water interface is used where porosity abruptly jumps from the soil porosity to unity (pure fluid). While this approach works well in fully explicit double-point MPM \citep{Bandara2015, Liu2017} and semi-implicit formulations without free water \citep{Xie2025}, the abrupt changes in porosity and permeability at the soil-water interface can cause divergence in the soil constitutive model when free water is present in a semi-implicit approach.

\begin{figure}[htbp]
	\begin{center}
		\includegraphics[width=0.95\textwidth]{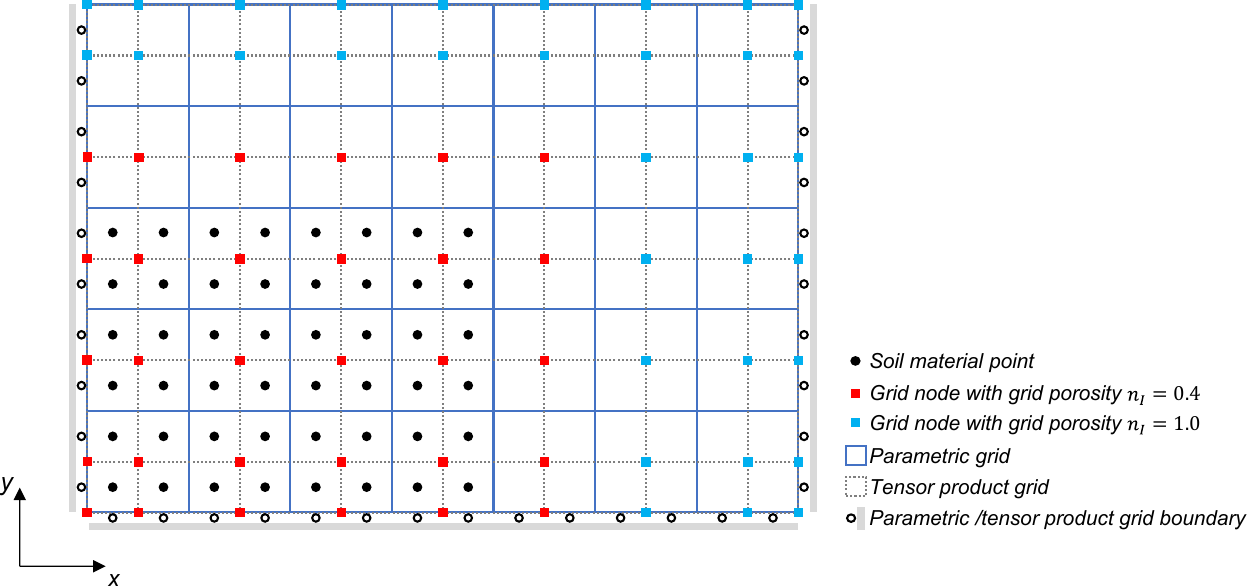}
		\caption{Illustration of a sharp soil-water interface}
		\label{DistinctInterface}
	\end{center}
\end{figure}

For a sharp interface, the grid porosity $\pmb{n}_I$ is obtained through mass-weighted mapping:
\begin{equation}\label{eq:nsp_nI}
	\pmb{n}_{I} = \frac{1}{\pmb{M}_{s I} } \sum_{sp=1}^{{N}_{sp}} {S}_{Isp} \, m_{s p} \, n_{sp} \,,
\end{equation}
where ${N}_{sp}$ is the total number of soil material points associated with grid node $I$; $m_{s p}$ is the soil material point mass; $n_{sp}$ is the soil material point porosity; and
\begin{equation}
	\pmb{M}_{s I} = \sum_{s p=1}^{{N}_{s p}} {S}_{Is p} \, m_{s p}
\end{equation}
is the lumped soil nodal mass. Nodes with zero soil mass are assigned unit porosity. The grid permeability $\pmb{k}_I$ is obtained similarly. This sharp interface, illustrated in Fig.~\ref{DistinctInterface}, introduces significant instability when water enters the soil domain due to suddenly applied drag forces.

\begin{figure}[htbp]
	\begin{center}
		\includegraphics[width=0.95\textwidth]{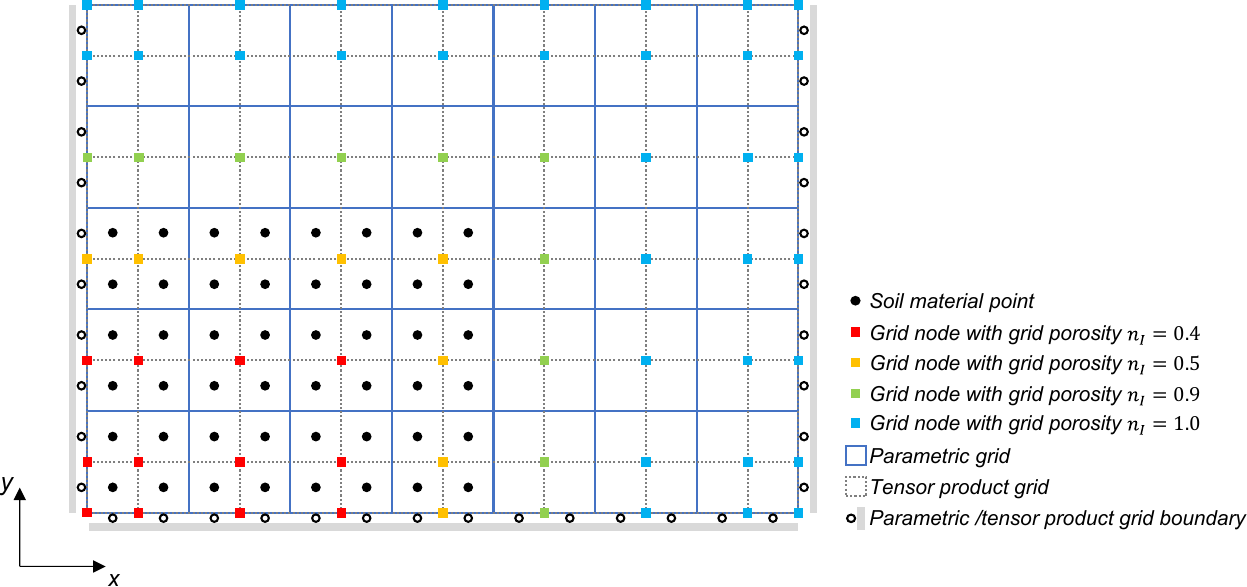}
		\caption{Illustration of a smooth soil-water interface}
		\label{smoothInterface}
	\end{center}
\end{figure}

To address this, \citet{Yamaguchi2020,Pan2021,Chandra2024b} introduced a smooth soil-water interface by mapping porosity as:
\begin{equation}\label{eq:nsp_nI_smooth}
	\pmb{n}_{I} = 1 - \frac{1}{V_{I} } \sum_{sp=1}^{{N}_{sp}} {S}_{Isp} \, V_{s p} \, (1 - n_{sp}) \,,
\end{equation}
where the grid nodal volume $V_{I}$ is obtained as
\begin{equation}\label{eq:nodal_volume}
    V_I = \sum_{gp=1}^{{N}_{gp}} {S}_{Igp}^0 \, V_{g p}^0 \;,
\end{equation}
through element-wise integration performed initially over the entire computational domain using Gaussian quadrature, with quadrature points matching the material point arrangement. The nodal volume $V_{I}$ is also utilised in the nodal-based free-water surface detection and $\delta$-correction methods (Sects.~\ref{sec:free_surface_detection} and~\ref{sec:delta_correction}). As illustrated in Fig.~\ref{smoothInterface}, this formulation creates transition layers where porosity gradually decreases from unity to the soil porosity.

While \citet{Yamaguchi2020,Pan2021,Chandra2024b} apply smooth porosity mapping but retain sharp permeability transitions, this research implements gradual permeability changes to further improve stability. This is achieved by evaluating the Kozeny-Carman equation \citep{Carrier2003} at grid nodes:
\begin{equation}
    \pmb{k}_{I} = C_{I}\frac{(\pmb{n}_{I})^3}{(1-\pmb{n}_{I})^2} \,,
\end{equation}
where
\begin{equation}
	C_{I} = \frac{1}{\pmb{M}_{s I} } \sum_{sp=1}^{{N}_{sp}} {S}_{Isp} \, m_{s p} \, C_{sp} \,,
\end{equation}
and 
\begin{equation}
	C_{sp} = k_{sp}\frac{(1-{n}_{sp})^2 }{({n}_{sp})^3}
\end{equation}
is a material constant obtained initially by inverting the Kozeny-Carman equation given the soil material point porosity $n_{sp}$ and permeability $k_{sp}$. For soils where permeability is not directly measured, $k_{sp}$ (cm/s) can be estimated from porosity $n_{sp}$ and effective particle diameter $D_{eff}$ (cm) at 20 °C \citep{Carrier2003}:
\begin{equation}\label{eq:est_permeability}
    k_{sp} = \frac{1.99\times10^4 \, (e^3)}{{(SF/D_{eff})}^2(1+e)} \,,
\end{equation}
where
{\begin{equation}\label{eq:n_to_e}
    e = -n_{sp}/(n_{sp} - 1)
\end{equation}}is the void ratio, $SF$ is the shape factor (spherical—6.0; rounded—6.6; medium angularity—7.5; angular—8.4) \citep{Carrier2003}, and $D_{eff}$ is calculated from the particle size distribution \citep{Carrier2003}:
\begin{equation}\label{eq:Deff}
    D_{eff} = \frac{100\%}{\sum(f_i/D_{ave\,\,i})} \,,
\end{equation}
where $f_i$ is the particle fraction between two sieve sizes and $D_{ave\,\,i} = D^{0.5}_{li} \times D^{0.5}_{si}$ is the average particle size \citep{Carrier2003}.

{It is worth clarifying the physical interpretation of the smooth porosity transition introduced in this section. The porosity at each soil material point ($n_{sp}$) is updated using the determinant of its deformation gradient, representing the actual porosity of the soil. Volume changes and dilatancy of the granular phase are captured at the material point level through the constitutive model (NorSand in the present study, Sect.~\ref{sec:tsunami1} and \ref{sec:tsunami2}), and therefore, no additional closure is required for the granular phase. The smooth transition illustrated in Fig.~\ref{smoothInterface} arises only at the grid nodes through the volume-weighted mapping of soil material points via the higher-order B-spline shape functions in Eq.~\eqref{eq:nsp_nI_smooth}, and is therefore a numerical conceptualisation rather than a physical diffused interface. The width of the transition layer is directly tied to the support of the B-spline shape function and decreases as the grid is refined.}

\subsection{A nodal-based free-water surface detection method}\label{sec:free_surface_detection}

In a semi-implicit coupled MPM, the water phase is implicit. Therefore, the Dirichlet boundary condition (zero water pressure) needs to be imposed on the free-water surface when obtaining the water pressure. Existing studies \citep{Yamaguchi2020,Chandra2024,Chandra2024b,He2024} detected the free-water surface using a proposed cell-based free-water detection method \citep{Kularathna2017}. In this method, the cells that are fully filled with water (that is, the water body) are determined by
\begin{equation}\label{eq:free_surface_cell}
    \frac{\sum_{wp=1}^{N_{wp}\in e} V_{wp}}{V_e} > \epsilon_\theta \,,
\end{equation}
where $\epsilon_\theta$ is a tolerance close to 1.0 (i.e. when a cell is full of water, this tolerance has to be 1.0 in theory). Zero water pressure can be applied to the nodes associated with cells that are not fully filled with water.

Although the cell-based free-surface detection method has been widely used, binary free-surface pressure imposition errors can be introduced, resulting in numerical instability \citep{Chandra2024b}. It is necessary to determine the appropriate tolerance to decide whether the cell is filled with water. This tolerance is set to 0.9 in \citet{Kularathna2017}, using the linear shape function. However, $\epsilon_\theta$ is reduced to 0.75 for a higher-order B-spline shape function in \citet{Chandra2024b} to ensure the stability of the free-water surface, as this type of shape function covers multiple cells, requiring a more lenient tolerance. However, 0.75 appears to be a compromised value derived from a trial-and-error process that deviates significantly from the ideal value of 1.0. Furthermore, in a knot-vector-based B-spline MPM, the cell's vertices do not overlap with the tensor product grid nodes (control points), as illustrated in Fig.~\ref{QuadraticBsplineGrids}. Therefore, a better algorithm is needed to determine the free-water surface, especially for MPM with a higher-order shape function.

\begin{figure}[htbp]
	\begin{center}
		\includegraphics[width=0.95\textwidth]{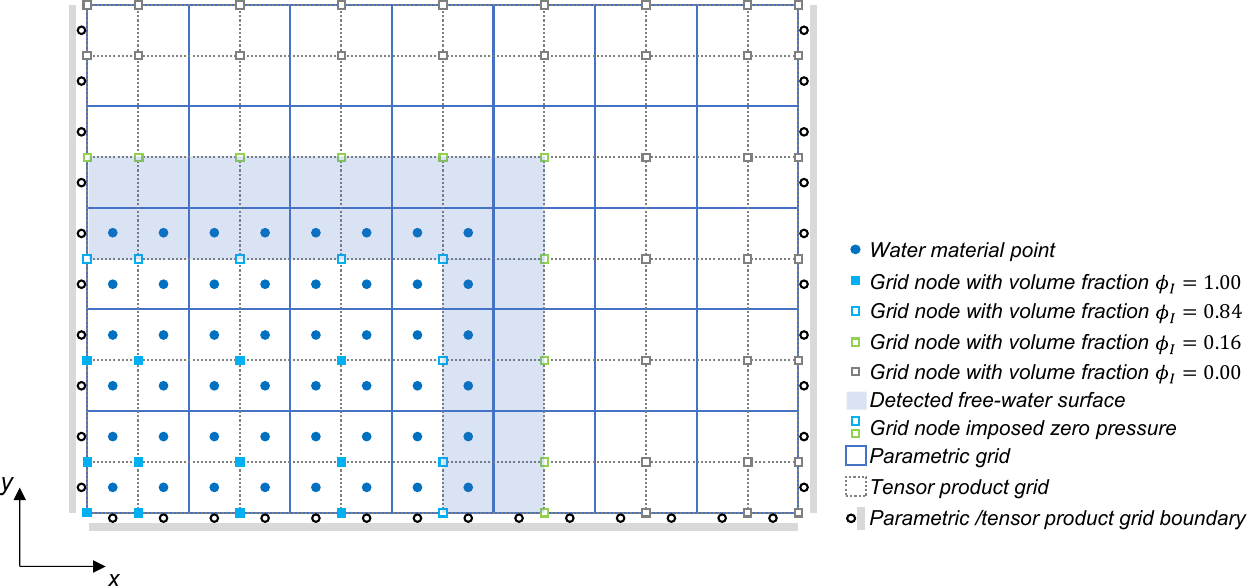}
		\caption{A schematic illustration of the nodal-based free-water surface detection method}
		\label{schematic2mp}
	\end{center}
\end{figure}

Since zero water pressure is applied on the grid nodes, the use of a nodal-based approach, rather than a cell-based algorithm, seems more reasonable. Therefore, in this research, a new nodal-based free-surface detection method has been developed, based on the ideas of creating a smooth soil-water interface and the cell-based method. In this nodal-based method, the grid nodes that are completely full of water are determined by calculating the volume fraction $\phi_I$ at each grid node:
\begin{equation}\label{eq:free_surface_node}
	\phi_I = \frac{\sum_{wp=1}^{{N}_{wp}} {S}_{Iwp} \, V_{w p}}{V_{I}} > \epsilon_\theta \,,
\end{equation}
where $V_{I}$ is the nodal volume previously defined in Eq.~(\ref{eq:nodal_volume}); $\epsilon_\theta$ is the tolerance similar to that used in the cell-based method. However, a stricter tolerance can be used for this nodal-based approach, resulting in more accurate and stable results. For example, a value of 0.99 is consistently used for the numerical examples presented hereinafter, which is significantly stricter than that used in the cell-based method as mentioned above. A schematic illustration of the proposed method is shown in Fig.~\ref{schematic2mp}. The blue-coloured grid nodes in Fig.~\ref{schematic2mp} are those detected as full of water (the water body), and zero water pressure is applied to the other nodes without a filled colour. The grid node's volume fraction gradually reduces from 1 (the water body) to 0 (void) at the free-water surface due to the loss of support of the shape function as the distance from material points to the outer nodes increases. Finally, the free-water surface at the grid-node level can be determined using the property of grid volume fraction as illustrated in Fig.~\ref{schematic2mp}.

\subsection{$\delta$-correction}\label{sec:delta_correction}

In a two-phase double-point MPM developed by \citep{Bandara2015} and \citep{Xie2025}, the volume of a water material point is updated by
\begin{equation}\label{eq:volume_update_wp}
    V_{wp}^{t + 1} = \left(n_{wp}^{0}/n_{wp}^{t + 1}\right) \, V_{wp}^{0} \,,
\end{equation}
where $V_{wp}$ and $n_{wp}$ are the volume and porosity of a water material point, respectively. The superscripts $t + 1$ and $0$ indicate these values in the current and initial time steps, respectively. When water flows through porous media, such as in a seepage flow problem, the volume of a water material point changes significantly when the water material point enters or leaves the porous medium. For example, consider a water material point, originally outside the porous medium with a porosity of 0.5, entering the porous medium. In this case, according to Eq.~(\ref{eq:volume_update_wp}), the volume of this material point will double immediately after entering, because the porosity of this material point at the beginning is 1.0 (free water) and it contracts to 0.5 after entering the porous medium. Due to the sudden increase in the volume of the water material point, these water material points entering the porous medium will overlap without spreading, resulting in a significant quadrature error.

\citet{Baumgarten2023} proposed the $\delta$-correction method to redistribute the clumped material points that cause the quadrature error. To minimise the global nodal volume error ${(\norm{E}^2)}_w$, this method incrementally corrects the coordinates of the water material points $\pmb{x}_{wp}$ by moving them along the gradient of the norm error with respect to each material point ${(\nabla\norm{E}^2)}_{wp}$ \citep{Baumgarten2023}. The global nodal volume error ${(\norm{E}^2)}_w$ can be obtained by adding the nodal volume error $E_{wI}$ in all active grid nodes:
\begin{equation}
    {(\norm{E}^2)}_w = \sum^{N_n}_{I = 1} E_{wI}^2 \,,
\end{equation}
where
\begin{equation}\label{eq:E_I}
    E_{wI} = \max(0, -V_I + \sum_{wp=1}^{{N}_{wp}} {S}_{Iwp} \, V_{w p}) \,.
\end{equation}
Equation~(\ref{eq:E_I}) is used to find how much the volume of current material points accumulated at each node exceeded the predefined nodal volume. As shown in Eq.~(\ref{eq:E_I}), in the case of non-exceeding, the nodal volume error, $E_I$ becomes zero.

Similarly, the gradient of the norm error with respect to each material point ${(\nabla\norm{E}^2)}_{wp}$ can be obtained as
\begin{equation}
    {(\nabla\norm{E}^2)}_{wp} = 2 V_{wp} \sum^{N_n}_{I = 1} \nabla{S}_{Iwp} \, E_{wI} \,.
\end{equation}

Then, the coordinates of the water material point are corrected by moving backwards a step along the gradient ${(\nabla\norm{E}^2)}_{wp}$:
\begin{equation}\label{eq:delta_correction}
    \pmb{x}_{wp}^{t + 1} \leftarrow \pmb{x}_{wp}^{t + 1} - b_0{(\nabla\norm{E}^2)}_{wp} \,,
\end{equation}
where
\begin{equation}
    b_0 = \frac{{(\norm{E}^2)}_w}{ \sum_{wp=1}^{{N}_{wp}} [{(\nabla\norm{E}^2)}_{wp} \cdot {(\nabla\norm{E}^2)}_{wp})] }
\end{equation}
is the step size.

In this research, we use the 1-step $\delta$-correction using Eq.~(\ref{eq:delta_correction}), which has been proven to be reliable and computationally efficient by \citet{Chandra2024,Chandra2024b}. Alternatively, a more accurate correction can be achieved using an iterative approach, as detailed in \citet{Baumgarten2023} (which was not applied in this research).

\section{Numerical examples}\label{sec:num_examples}
In this section, some numerical examples are presented to validate the 
newly proposed formulations for modelling soil-water coupled problems 
with free-surface flows. Validations are conducted against several 
experiments from the literature, starting from a classic free-water 
dam break problem reported by \citet{Lobovsky2014}. Then, the 
performance of the model for seepage flow has been validated with the 
experiments conducted by \citet{Liu1999}. Two granular flow 
triggered tsunami experiments, conducted by \citet{Rauter2022} and 
\citet{Sarlin2022}, respectively, are simulated using the proposed 
method with NorSand constitutive model, using model parameters 
calibrated from triaxial tests from the literature. {Finally, a dam-break wave over a movable granular bed \citep{Spinewine2005} is simulated to demonstrate the reverse interaction: fluid-induced granular deformation.} {It is worth noting that these examples also showcase advantages of the proposed scheme over traditional fractional step approaches: the incremental fractional step method preserves the water volume across the relatively long simulation times required, removes the permeability-dependent time-step restriction (which becomes prohibitively expensive for the lower permeabilities typical of real geotechnical materials, as discussed in \citet{Xie2025}), and provides the numerical stability needed to integrate advanced constitutive models such as NorSand in the last three examples.}

\subsection{Dam break}\label{sec:dambreak}

\begin{figure}[htbp]
	\begin{center}
		\includegraphics[width=0.45\textwidth]{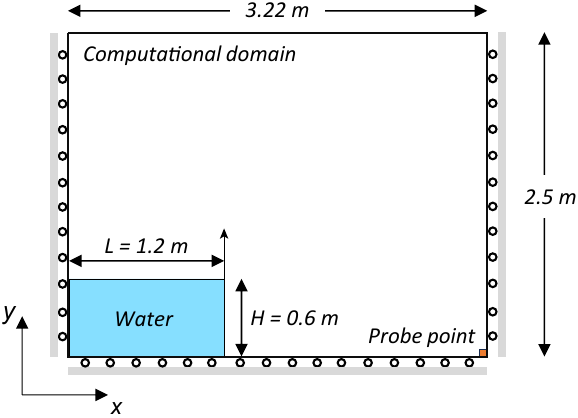}
		\caption{Dam break flow - setup of the initial geometry and boundary conditions}
		\label{dambreak_setup}
	\end{center}
\end{figure}

A classic dam break flow over a dry bed problem is used to investigate the performance of the proposed method in simulating free-water flow. In the dam break experiment conducted by \citet{Lobovsky2014}, a water column is held by a removable gate in a flume, as illustrated in Fig.~\ref{dambreak_setup}. After the sudden removal of the gate, the water column collapses, and water flow is generated. After the water flow reaches the end of the flume, a plunging breaking wave is generated. The water pressure time history was recorded at the probe point (Fig.~\ref{dambreak_setup}) during this process, which is compared with the results obtained with the numerical simulations using the proposed methodology.

The numerical model is set up based on the geometry and boundary conditions shown in Fig.~\ref{dambreak_setup}. The computational domain (the background mesh) has dimensions of 3.22~m in length and 2.5~m in height. A plane strain condition is adopted because the width (thickness) of the flume is sufficiently large and the inward deformation is negligible. Roller boundary conditions are applied on the sides and bottom of the computational domain to represent a perfectly smooth flume wall. The water column's length $L$ and height $H$ are 1.2 and 0.6~m, respectively. A structured quadrilateral background mesh with a cell size of 0.02~m is used to discretise the computational domain, and $4^2$ water material points are uniformly placed in each cell at the initial configuration to discretise the water body, resulting in 28,800 water material points in total. The removable gate is not modelled for this numerical example since it was removed rapidly at time $t$ = 0~s during the experiment \citep{Lobovsky2014}. Instead, the gravitational acceleration $g$ = 9.81~$\text{m/s}^2$ is applied to the numerical model at $t$ = 0~s to generate the dam break water flow. A time increment $\Delta t$ = 0.0001~s (about 0.99~CFL) is used for this dam break problem, and the final simulation time is taken as $t$ = 7.00~s.

In this dam break flow problem, we investigate the performance of the proposed numerical method by carrying out simulations using the FLIP and TPIC (1\%) schemes. The performance of the hyperelastic water model has also been explored in this problem. {
The water stress is updated using Eq.~\eqref{eq:sig_w}. For the simulation without the hyperelastic water model, the water contribution is represented directly by the classical hydrostatic term \(p_w\boldsymbol{I}\), as in the baseline formulation.
}

\begin{figure}[htbp]
	\begin{center}
		\includegraphics[width=0.95\textwidth]{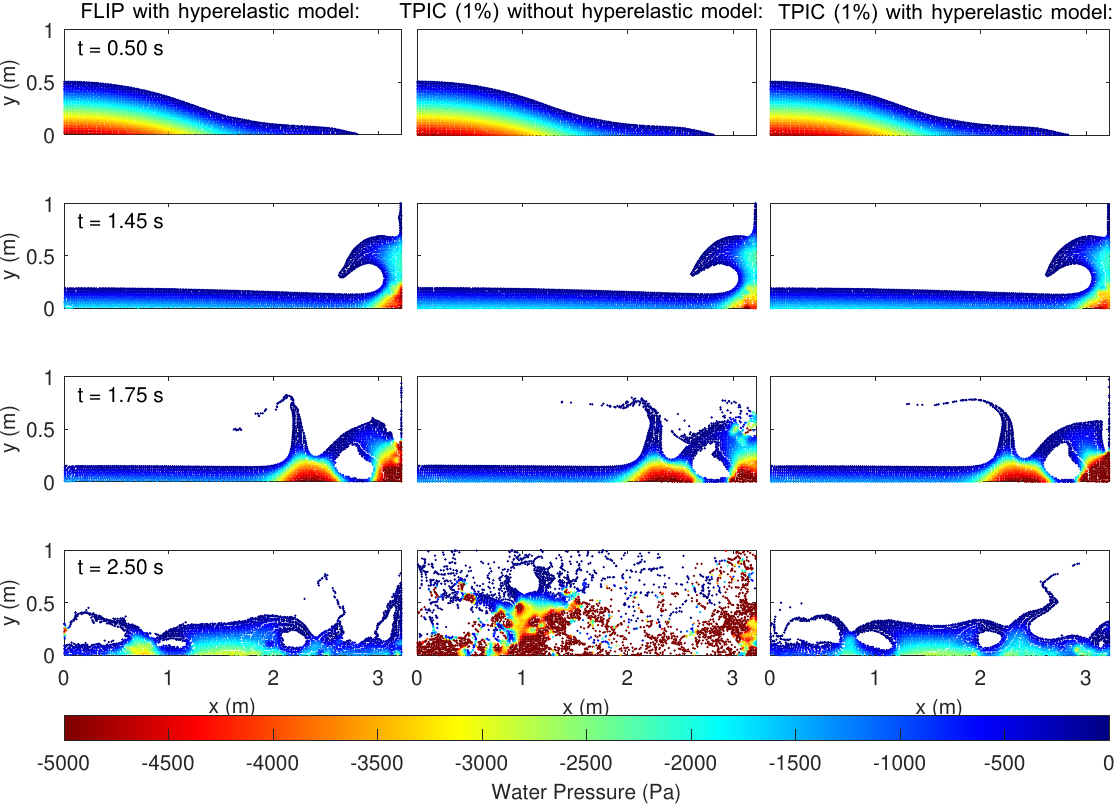}
		\caption{Dam break flow - the evolution of pressure field at different time snapshots t = 0.50, 1.45, 1.75, and 2.50~s (non-dimensional time $t\sqrt{g/H}$ = 2.02, 5.86, 7.08, and 10.11)}
		\label{dambreak_compare}
	\end{center}
\end{figure}

Figure~\ref{dambreak_compare} presents the evolution of the dam break flows at different time snapshots. Note that, for the dam break problem, a non-dimensional time $t\sqrt{g/H}$ is usually used to ensure a normalised result valid over various scales of experiments. The simulations are given by the FLIP scheme with the hyperelastic water model, the TPIC (1\%) scheme without the hyperelastic water model, and the TPIC (1\%) scheme with the hyperelastic water model. Before the water impacts the right end of the flume, all three methods show identical results. At $t$ = 1.45~s, the plunging breaking wave has been generated, and the waves given by FLIP and TPIC (1\%) schemes have slightly different geometries at the wave's lip. At $t$ = 1.75~s, the first wave falls, and the second wave is generated. Some severe instabilities can be observed in the simulation without the hyperelastic water model. However, the hyperelastic water model does not change the dynamics of the water, maintaining the same shape as the one without it, whilst stabilising the numerical model. The TPIC (1\%) shows more stable results than the FLIP scheme. At $t$ = 2.50~s, the simulation given by TPIC (1\%) without the hyperelastic water model has completely diverged due to strong instabilities. In contrast, both FLIP and TPIC (1\%) schemes with hyperelastic water models can produce very stable solutions.

\begin{figure}[htbp]
	\begin{center}
		\includegraphics[width=0.95\textwidth]{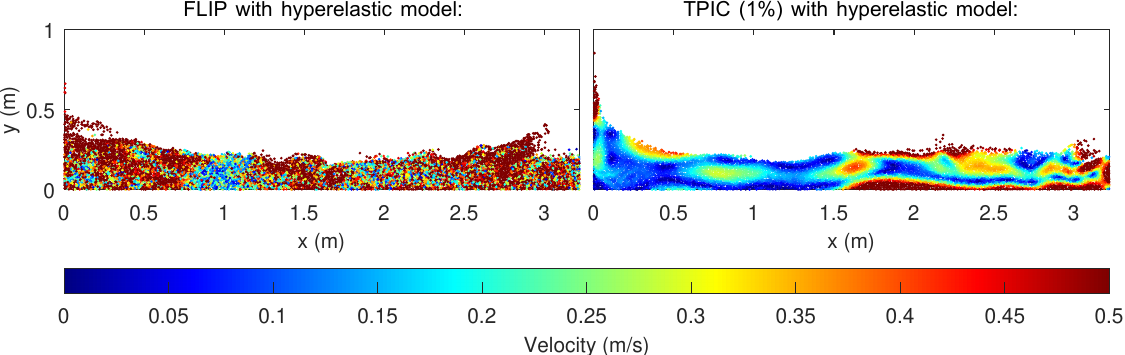}
		\caption{Dam break flow - the velocity fields given by the FLIP and TPIC (1\%) schemes with hyperelastic water model at t = 7.00~s (non-dimensional time $t\sqrt{g/H}$ = 28.30)}
		\label{dambreak_compare_velocity}
	\end{center}
\end{figure}

The comparison of the velocity field of the FLIP and TPIC (1\%) schemes at 7.00~s is shown in Fig.~\ref{dambreak_compare_velocity}. The velocity field of the FLIP scheme shows very strong oscillations despite its displacement field being relatively stable. The pattern of this oscillation is similar to volumetric locking, which occurred in the stress field. As shown in Fig.~\ref{dambreak_compare_velocity}, the TPIC (1\%) scheme is sufficient to eliminate the velocity oscillation in MPM, resulting in a smooth contour map.

\begin{figure}[htbp]
	\begin{center}
		\includegraphics[width=0.95\textwidth]{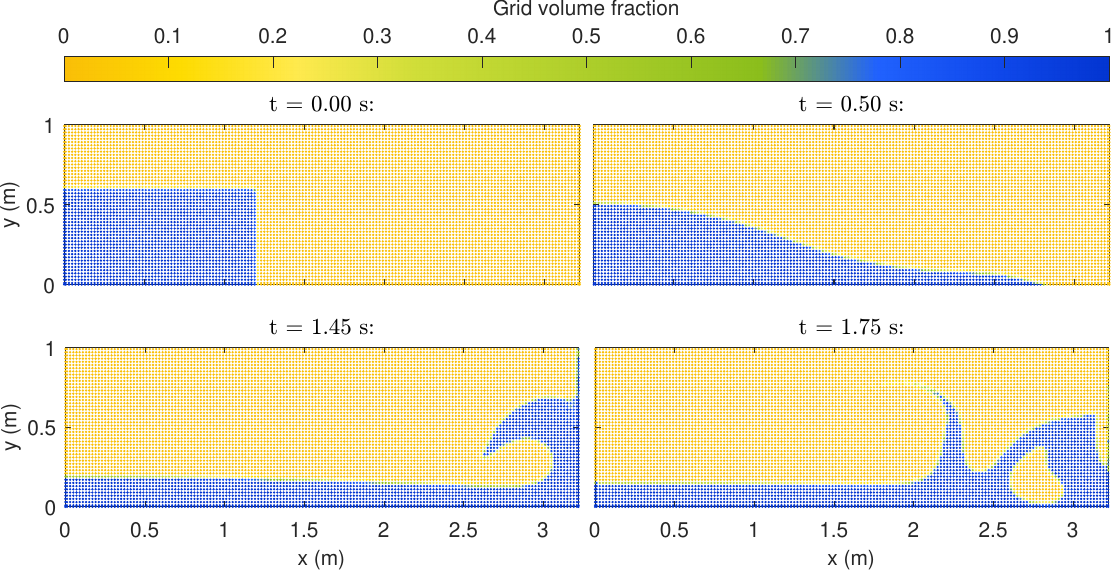}
		\caption{Dam break flow - the evolution of grid volume fraction for TPIC (1\%) schemes with hyperelastic water model at different time snapshots t = 0.00, 0.50, 1.45, and 1.75~s (non-dimensional time $t\sqrt{g/H}$ = 0.00, 2.02, 5.86, and 7.08)}
		\label{dambreak_gridVfraction}
	\end{center}
\end{figure}

The excellent performance of the proposed nodal-based free-water surface detection method is demonstrated by Fig.~\ref{dambreak_gridVfraction}. The complex wave geometry has been tracked precisely using the proposed method, ensuring accurate imposition of zero water pressure.

\begin{figure}[htbp]
	\begin{center}
		\includegraphics[width=0.95\textwidth]{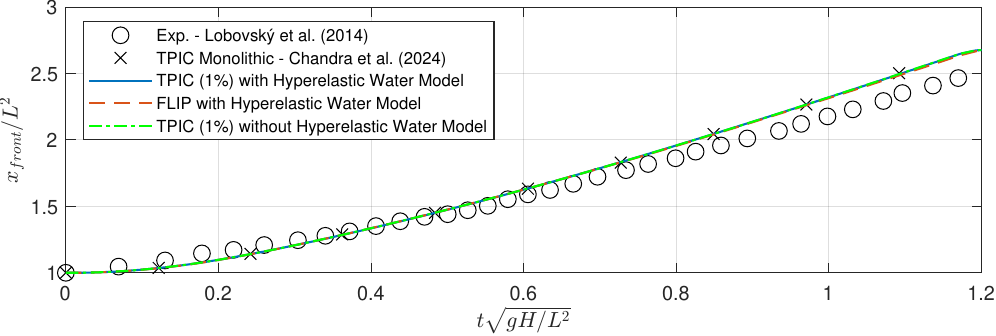}
		\caption{Dam break flow - time evolution of the waterfront position: comparisons of the experiment conducted by \citet{Lobovsky2014} and the numerical simulations given by \citet{Chandra2024}, our proposed method with FLIP and TPIC (1\%) schemes}
		\label{dambreak_compare_runout}
	\end{center}
\end{figure}

\begin{figure}[htbp]
	\begin{center}
		\includegraphics[width=0.95\textwidth]{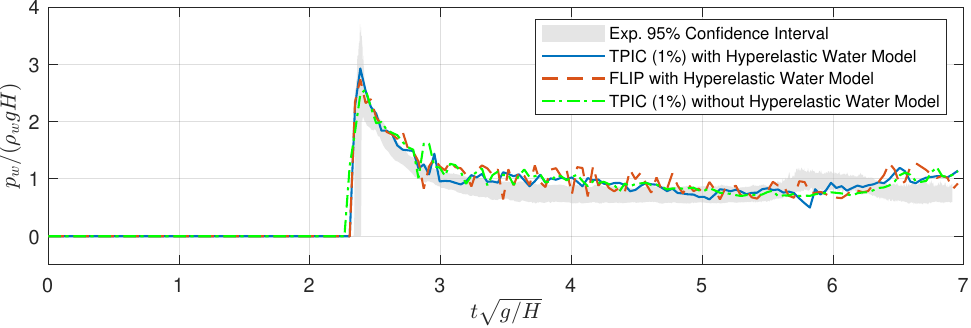}
		\caption{Dam break flow - the recorded pressure time history at the probe point: comparisons of the experiment conducted by \citet{Lobovsky2014} and the numerical simulations given by FLIP and TPIC (1\%) schemes}
		\label{dambreak_compare_pressure}
	\end{center}
\end{figure}

Figure~\ref{dambreak_compare_runout} compares the time evolution of the waterfront positions given by the experiment, the previous numerical study using monolithic MPM with TPIC scheme conducted by \citet{Chandra2024}, and with the new method proposed in this research. As shown in Fig.~\ref{dambreak_compare_runout}, the proposed method yields identical results to the previous numerical study and shows excellent agreement with the experiment. Importantly, the proposed method achieves this level of accuracy with significantly lower computational cost than the Monolithic approach of \citet{Chandra2024}, whilst maintaining the same numerical stability. Figure~\ref{dambreak_compare_pressure} compares the recorded water pressure time history at the probe point. The water pressure time history recorded by the transducer during the experiment has relatively large uncertainty compared to the evolution of the waterfront positions. \citet{Lobovsky2014} repeated the experiment 100 times, and a 95\% confidence interval for those 100 repetitions is provided. As shown in Fig.~\ref{dambreak_compare_pressure}, our proposed method also demonstrates excellent agreement with the experiment regarding the water pressure during the formation of complex breaking waves. Specifically, the TPIC (1\%) scheme shows less noise in the temporal distribution of the water pressure compared to the FLIP scheme. Additionally, the hyperelastic water model slightly reduces the temporal pressure vibration with a negligible increase in computational cost.

\clearpage

\subsection{Dam break seepage flow}\label{sec:seepage}

\begin{figure}[htbp]
	\begin{center}
		\includegraphics[width=0.6\textwidth]{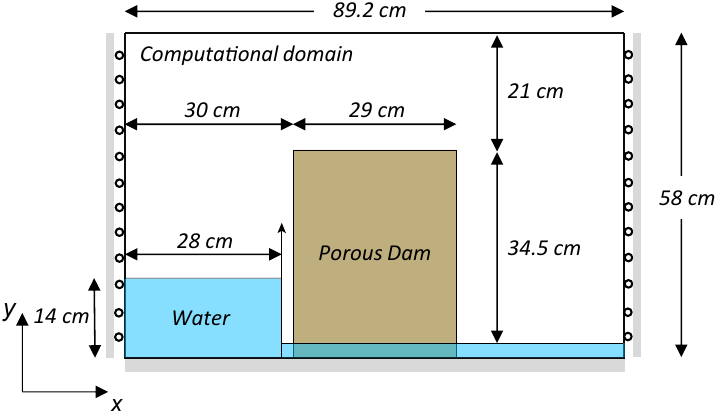}
		\caption{Dam break seepage flow - setup of the initial geometry and boundary conditions}
		\label{dambreak_seepage_setup}
	\end{center}
\end{figure}

The second numerical example validates the proposed numerical method against a dam break seepage flow experiment reported by \citet{Liu1999}. The experiment setup is illustrated in Fig.~\ref{dambreak_seepage_setup}. The flume used in the experiment conducted by \citet{Liu1999} has dimensions of 89.2~cm in length and 58~cm in height. A 29~cm high and 37~cm long porous dam is located 30~cm from the left end of the flume. The porous dam is formed of glass beads with a 0.3~cm uniform diameter and a porosity of 0.39. These glass beads are bonded with adhesive so they can remain fixed in place during the test \citep{Liu1999}. A 2 cm thick gate is placed on the left-hand side of the porous dam to form an upstream reservoir with a 14 cm water depth. The downstream water depth is 2.5~cm in the initial setup. According to \citet{Liu1999}, the gate was lifted manually within 0.1~s at $t$ = 0~s, and a camera with a frame rate of 10 frames per second (fps) is used to record the upcoming seepage flow.

The geometry and boundary conditions for the numerical model closely reproduce the conditions of the experiment. The porous dam is modelled as an elastic material fixed at the bottom of the flume with the same porosity as the porous medium in the experiment. Given the initial porosity $n_0$ and the diameter of these glass beads, the initial permeability $k_0$ for this porous dam can be calculated using Eq.~(\ref{eq:est_permeability}), which is 0.0793 m/s. Similar to the previous numerical example, the release of the gate is not modelled, and the gravitational acceleration $g$ = 9.81~$\text{m/s}^2$ is applied to the numerical model at t = 0 s to simulate its quick removal. A structured quadrilateral background mesh with a cell size of 0.01 m discretises the computational domain; $2^2$ and $4^2$ material points are uniformly placed in each cell at the initial configuration to discretise the soil and water bodies, respectively. The number of soil material points per cell is significantly smaller than the number of water material points because the porous dam is rigid in place, so that the soil material points always remain in their optimised positions without producing quadrature error. Therefore, using a large number of soil material points per cell is not needed. A time increment $\Delta t$ = $5\times10^{-5}$~s (about 0.99~CFL) is used for this dam break problem, and the final simulation time is 4.0 s.

\begin{figure}
	\begin{center}
		\includegraphics[width=0.95\textwidth]{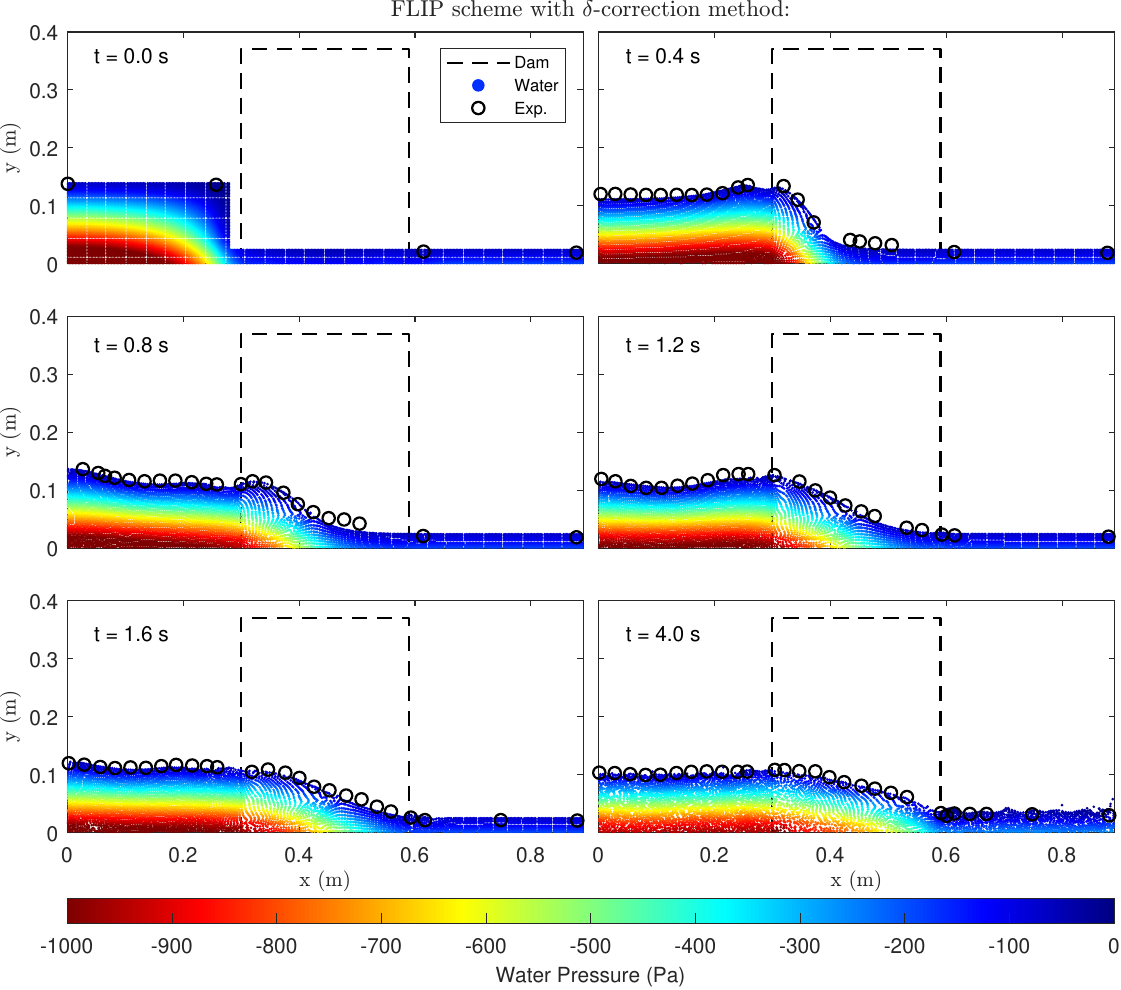}
		\caption{Dam break seepage flow - the evolution of pressure field at different time snapshots t = 0.0, 0.4, 0.8, 1.2, 1.6 and 4.0 s: results are given by FLIP scheme with $\delta$-correction method; "Exp." in the legend denotes the experimental measurements adopted from \citet{Liu1999}}
		\label{dambreak_seepage_FLIP_WaterModel}
	\end{center}
\end{figure}

\begin{figure}
	\begin{center}
		\includegraphics[width=0.95\textwidth]{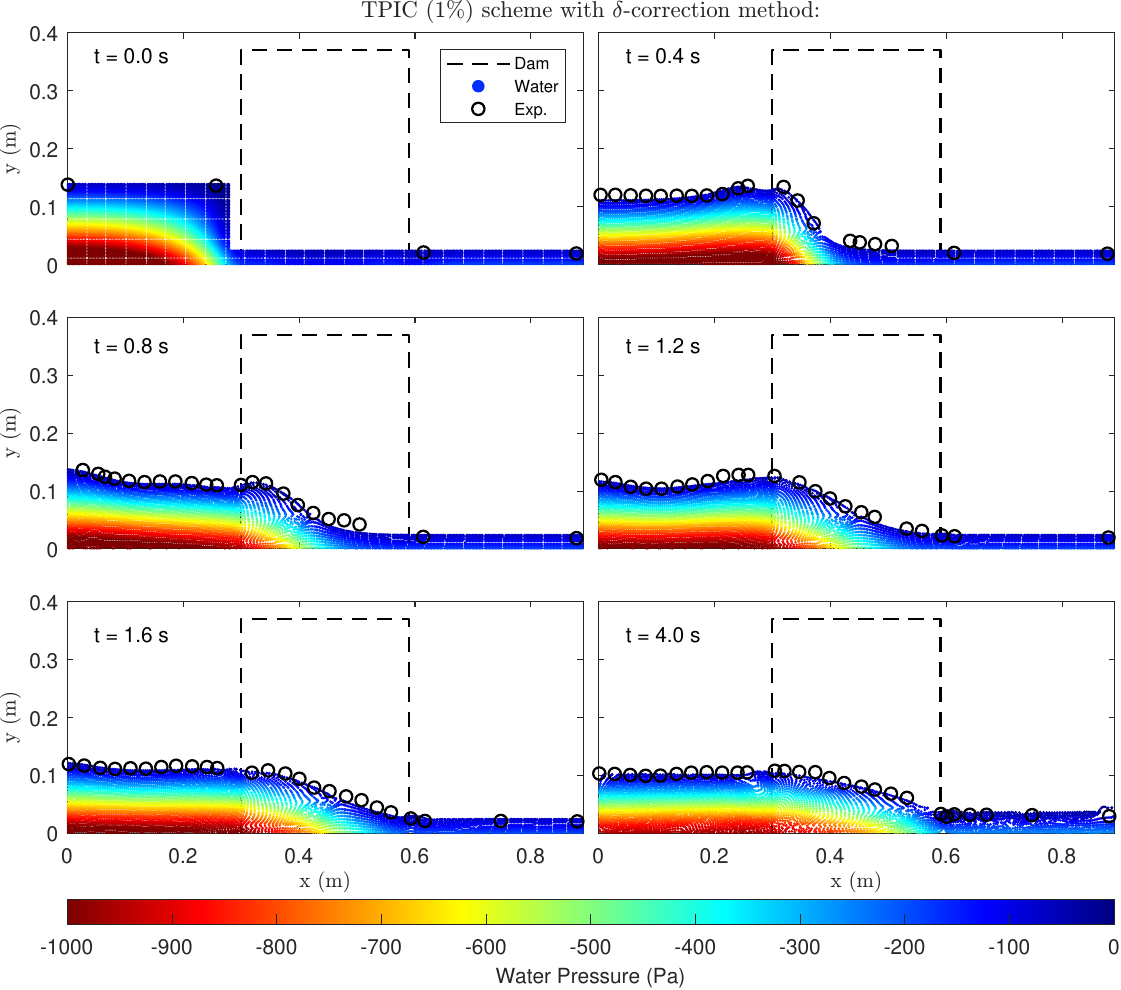}
		\caption{Dam break seepage flow - the evolution of pressure field at different time snapshots t = 0.0, 0.4, 0.8, 1.2, 1.6 and 4.0 s: results are given by TPIC (1\%) scheme with $\delta$-correction method; "Exp." in the legend denotes the experimental measurements adopted from \citet{Liu1999}}
		\label{dambreak_seepage_1pTPIC_WaterModel}
	\end{center}
\end{figure}

\begin{figure}
	\begin{center}
		\includegraphics[width=0.95\textwidth]{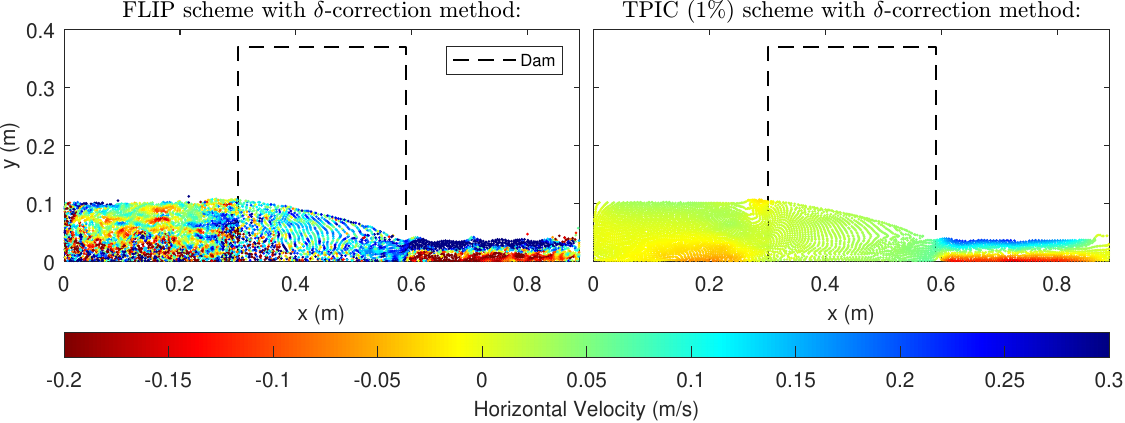}
		\caption{Dam break seepage flow - the horizontal velocity fields given by the FLIP and TPIC (1\%) schemes at t = 4.0 s}
		\label{dambreak_seepage_FLIPvs1pTPIC}
	\end{center}
\end{figure}

In this numerical example, the hyperelastic water model is consistently used for all simulations to ensure numerical stability. Figures~\ref{dambreak_seepage_FLIP_WaterModel} and \ref{dambreak_seepage_1pTPIC_WaterModel} present the comparison between the experiment and the numerical simulations given by FLIP and TPIC (1\%) schemes with $\delta$-correction, respectively. As shown in these figures, both FLIP and TPIC (1\%) schemes with $\delta$-correction method match the experiment very well, indicating the correctness of the proposed numerical model and the approach to estimate the permeability. The TPIC (1\%) scheme demonstrates a more stable but very similar water surface compared to that given by the FLIP scheme at $t$ = 4.0 s. Similarly to the previous dam break example, with water only, the velocity field given by the TPIC (1\%) scheme is very smooth, whilst the one given by the FLIP scheme shows strong oscillation, as shown in Fig.~\ref{dambreak_seepage_FLIPvs1pTPIC}. The velocity oscillation does not prevent the FLIP scheme with $\delta$-correction from producing a solution comparable with the experiment in terms of water levels. Compared to the previous dam break problem, where a significant difference develops, the flow speed in this case is lower, and the breaking wave does not exist.

\begin{figure}
	\begin{center}
		\includegraphics[width=0.95\textwidth]{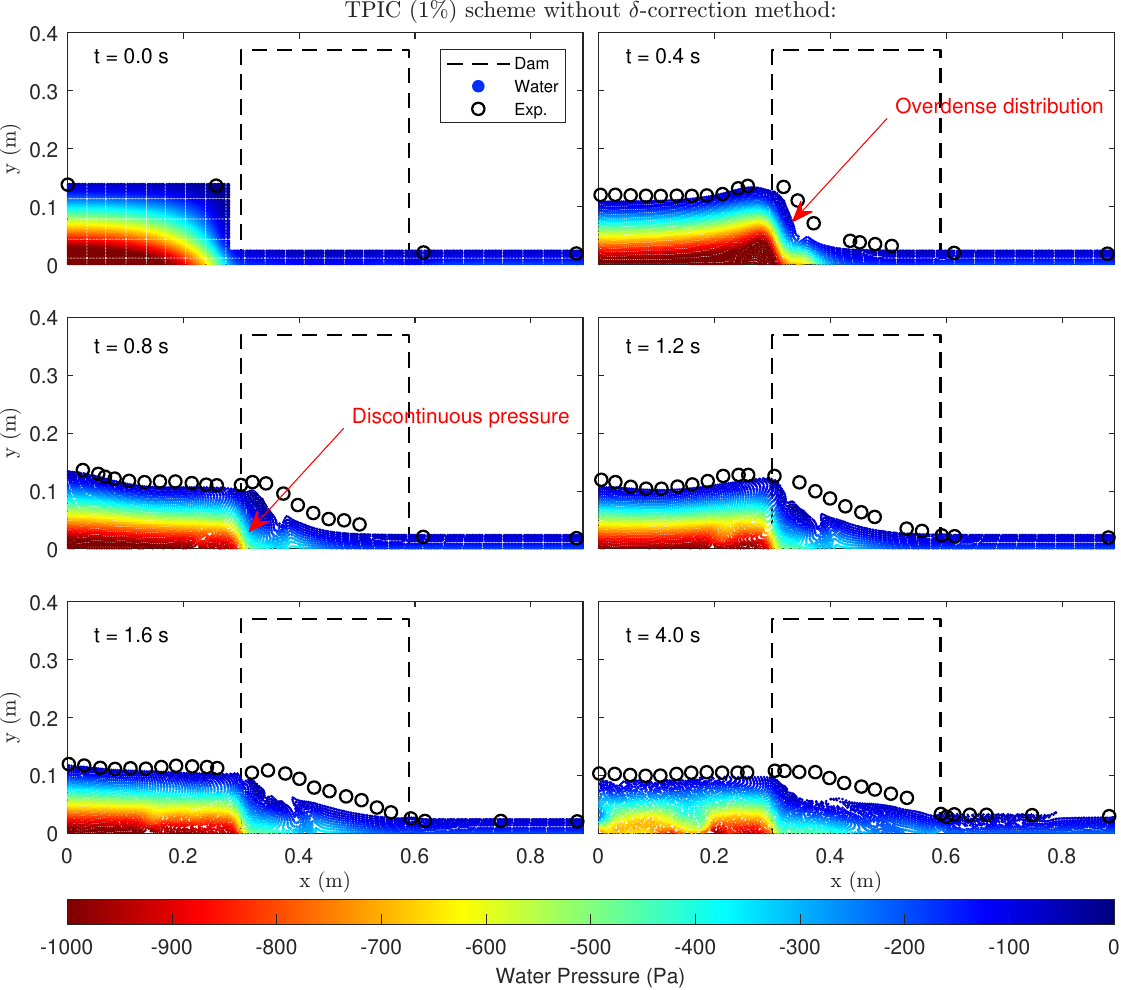}
		\caption{Dam break seepage flow - the evolution of pressure field at different time snapshots t = 0.0, 0.4, 0.8, 1.2, 1.6 and 4.0 s: results are given by TPIC (1\%) scheme without $\delta$-correction method; "Exp." in the legend denotes the experimental measurements adopted from \citet{Liu1999}}
		\label{dambreak_seepage_1pTPIC_WaterModel_noCorrection}
	\end{center}
\end{figure}

The performance of the $\delta$-correction method has been studied in this numerical example. Figure~\ref{dambreak_seepage_1pTPIC_WaterModel_noCorrection} compares the experiment and the results given by TPIC (1\%) schemes without the $\delta$-correction method. As we can observe, without the $\delta$-correction method, after the water enters the porous dam, the water material points do not have time to distribute properly because of the sudden enlargement of material point volume, resulting in an overcrowded distribution and discontinuous water pressure. Clearly, the TPIC (1\%) without the $\delta$-correction method cannot reasonably compare with the experiment, especially after the water flows into the porous medium. Therefore, implementing the $\delta$-correction method is necessary for problems where water flows rapidly in and out of porous media, such as the fast-moving seepage flow problem. 

{It is worth noting that, after redistribution by the 
$\delta$-correction, the spacing between water material points inside the porous dam is directly related to the local porosity. As described in Sect.~\ref{sec:delta_correction}, when a water material point enters a porous medium, its volume expands by a factor of $1/n$ relative to its free-water value, since the same mass of water now occupies only the pore space. This volume expansion translates into a sparser arrangement of water material points 
inside the porous domain, with the spacing scaling approximately as $(1/n)^{1/dof}$, where $dof$ is the spatial dimension. For the porous dam in this test with porosity 0.39 in 2D, the water material points inside the dam are expected to be approximately $1/\sqrt{0.39} \approx 1.60$ times the spacing of those in the free water region, consistent with the distribution shown in Figs.~\ref{dambreak_seepage_FLIP_WaterModel} and \ref{dambreak_seepage_1pTPIC_WaterModel}. This sparser distribution correctly reflects the fact that water fills only the pore space rather than the full domain.}

\clearpage

\subsection{Granular flow and tsunami: case 1}\label{sec:tsunami1}
\begin{figure}[htbp]
	\begin{center}
		\includegraphics[width=0.95\textwidth]{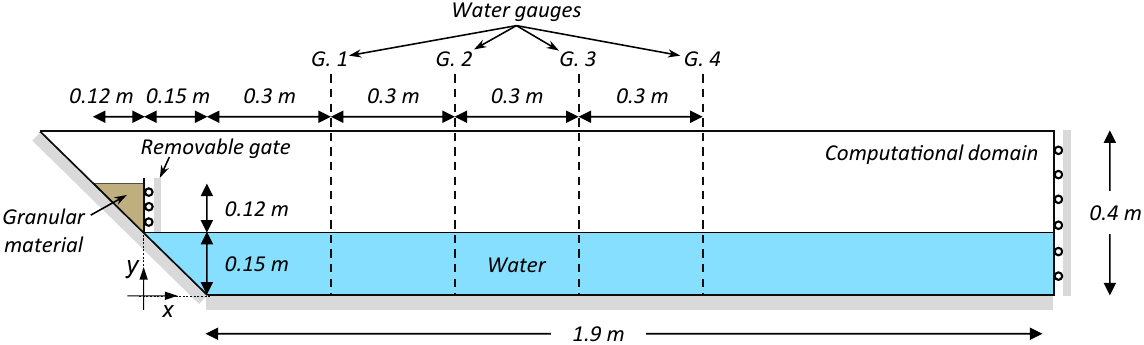}
		\caption{Granular flow and tsunami: case 1 - setup of the initial geometry and boundary conditions}
		\label{granular_flow_sunami_setup}
	\end{center}
\end{figure}

 In this numerical example, the proposed numerical model's performance, in terms of intense soil-water interaction, is studied. For this example, the TPIC (1\%) scheme with hyperelastic water model and $\delta$-correction method has been used. The proposed numerical method has been compared with the experiment conducted by \citet{Rauter2022}. In this experiment, granular material (glass beads with a 4 mm diameter) is held by a removable gate in an inclined slope inside a flume, which contains water to a depth of 0.15 m. This experimental setup is illustrated in Fig.~\ref{granular_flow_sunami_setup}. After suddenly removing the gate, the glass beads collapse and flow into the water, representing a landslide scenario that generates a tsunami. The wave height was recorded by four water gauges located in the flume with a 0.3 m interval, as shown in Fig.~\ref{granular_flow_sunami_setup}.

\begin{figure}
	\begin{center}
		\includegraphics[width=0.95\textwidth]{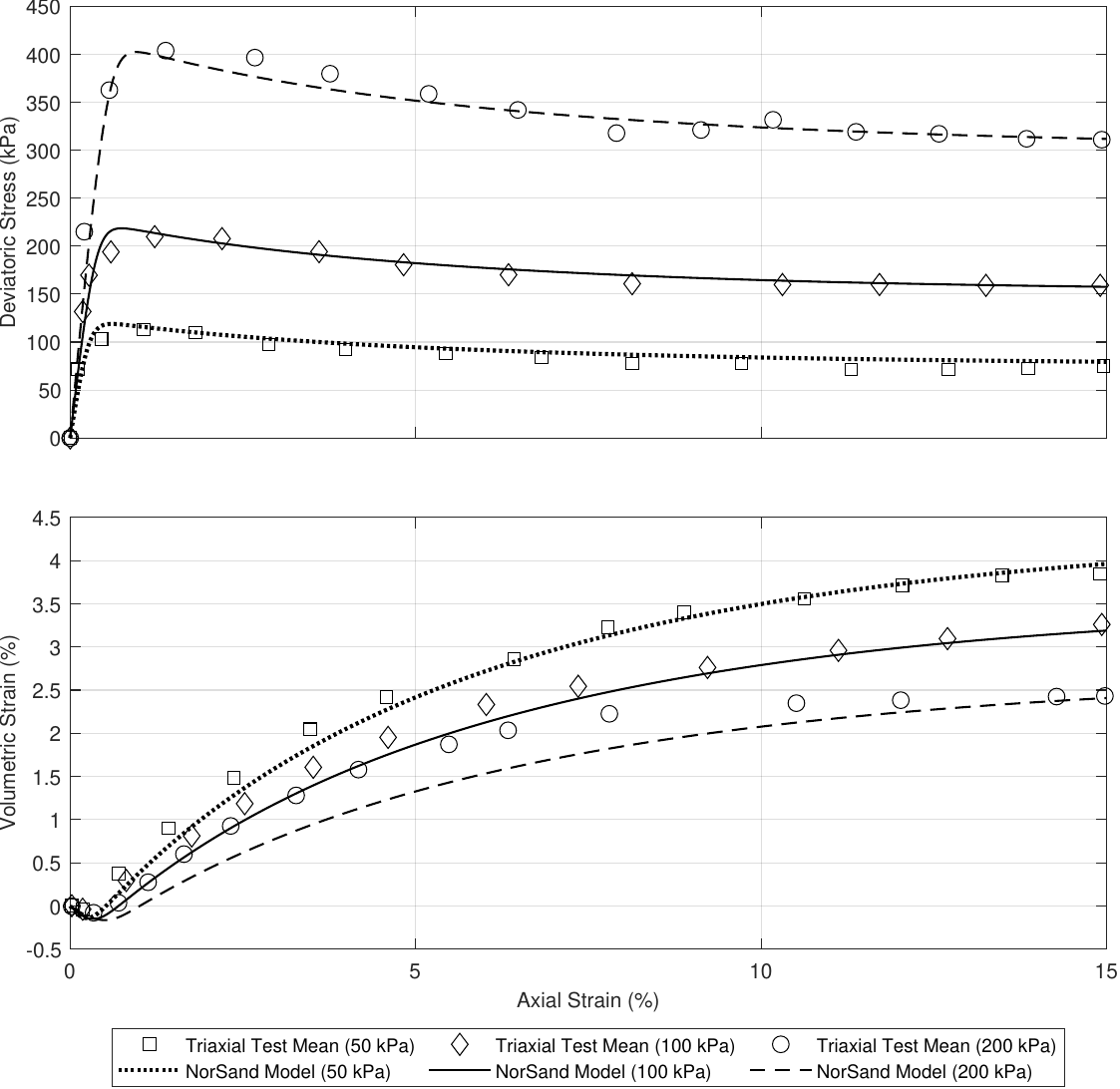}
		\caption{Granular flow and tsunami - calibration of the NorSand parameters for 4 mm diameter glass beads against the triaxial tests conducted by \citet{Wu2017}}
		\label{4mmglassbeads_calibration}
	\end{center}
\end{figure}

In this case, an advanced constitutive model, the hyperelastic NorSand model based on a large-strain framework \citep{Borja2006}, is adopted to properly simulate the complex behaviour of the granular material. The modified F-bar method developed in \citet{Xie2025} has been used to stabilise the solutions in the soil phase. \citet{Wu2017} conducted a series of drained triaxial tests for 4~mm diameter glass beads, which have been used in the research to calibrate the NorSand model parameters. In the drained triaxial tests conducted by \citet{Wu2017}, the samples were sheared under 50, 100 and 200~kPa confining pressures with a 0.545 initial void ratio. The tests were repeated three times for each confining pressure. The relationships between deviatoric stress, volumetric strain and axial strain were recorded \citep{Wu2017}, and we took the mean values of three tests for each confining stress to conduct this calibration. According to \citet{Wu2017}, the tested friction angle at the critical state $\phi_{cs}$ of those 4~mm glass beads is approximately 25.5°. Therefore, the slope of the critical state line can be calculated by $M = 6\sin{\phi_{cs}}/(3-\sin{\phi_{cs}})$ \citep{Jefferies2015}. The other model parameters are calibrated by fitting the curves. Table~\ref{tab:4mmglassbeads} presents the calibrated NorSand parameters after the trial and error procedure, and Fig.~\ref{4mmglassbeads_calibration} compares the triaxial tests given by the experiment and NorSand constitutive model. As shown in Fig.~\ref{4mmglassbeads_calibration}, the calibrated NorSand model can reproduce the experimental triaxial tests very accurately. Therefore, these calibrated parameters are used to simulate the landslide-tsunami experiment conducted by \citet{Rauter2022}, which uses the same granular material (i.e. 4~mm diameter glass beads). \citet{Rauter2022} reported that the grain density for those 4 mm glass beads is 2500~kg/$\text{m}^3$, and the initial packing fraction is approximately 0.6, which is equivalent to 0.4 porosity. Therefore, the initial void ratio $e_0$ can be obtained with Eq.~(\ref{eq:n_to_e}), yielding a value of 0.667. Then, the initial specific volume $v_0 = 1 + e_0 = 1.667$. Knowing the diameter and void ratio of the glass beads, the initial permeability $k_0$ can be calculated by Eq.~(\ref{eq:est_permeability}), which equals 0.157~m/s.

\begin{table}[htbp]
        \centering
	\caption{Granular flow and tsunami - NorSand model parameters for 4 mm glass beads and input parameters for numerical simulation}\label{tab:4mmglassbeads}%
	\begin{tabular}{@{}ll@{}}
		\toprule
		Parameter &Value  \\
		\midrule
		Water density $\rho_w \,[\text{kg/m}^3]$    & 1000 \\
		Soil grain density $\rho_s \,[\text{kg/m}^3]$    & 2500 \\
		Shear modulus $G \,[\text{MPa}]$    & 40 \\
		Swelling index $\kappa$    & 0.002 \\
            Reference specific volume $v_{c0}$    & 1.84 \\
            Compression index $\lambda$    & 0.02 \\
		Slope of the critical state line $M$    & 1.0 \\
		Yield function constant $N$    & 0.5 \\
            Plastic potential constant $\Bar{N}$    & 0.5 \\
		Hardening coefficient $h$    & 1000 \\
            Multiplier for maximum plastic dilatancy $\alpha$    & -11 \\
            Initial specific volume $v_0$    & 1.667 \\
		Initial porosity $n_0$    & 0.4 \\
		Initial hydraulic conductivity $k_0 \,[\text{m/s}]$    & 0.157 \\
		Gravitational acceleration $g \,[\text{m/s}^2]$    & -9.81 \\
		\hline
	\end{tabular}
\end{table}

The geometry and boundary conditions of the numerical model follow the setup shown in Fig.~\ref{granular_flow_sunami_setup}. The inclined slope is modelled as a fixed boundary condition. According to \citet{Rauter2022}, the gate is released by a 20~kg mass connected to the gate system, and the time $t$ = 0~s is when the gate fully opens. Therefore, in this case, it is necessary to model the process of releasing the gate. We model the gate as a roller boundary condition, and its length reduces with gravitational acceleration to represent the releasing process of the gate, which was lifted by a free-falling mass. The NorSand constitutive model must be initialised before the release of the gate. To do that, we apply the gravitational load gradually during 1~s to initialise the model, while the gate remains closed. For this numerical example, a structured quadrilateral background mesh with a cell size of 0.005~m is used to discretise the computational domain, and $4^2$ soil and water material points are uniformly placed in each cell according to the setup. This configuration is sufficient to ensure the accuracy of this numerical example. A time increment $\Delta t$ = $3\times10^{-5}$~s (about 0.95~CFL) is used, and the total simulation time is 3~s, excluding the time for initialisation and gate opening.

\begin{figure}
	\begin{center}
		\includegraphics[width=0.95\textwidth]{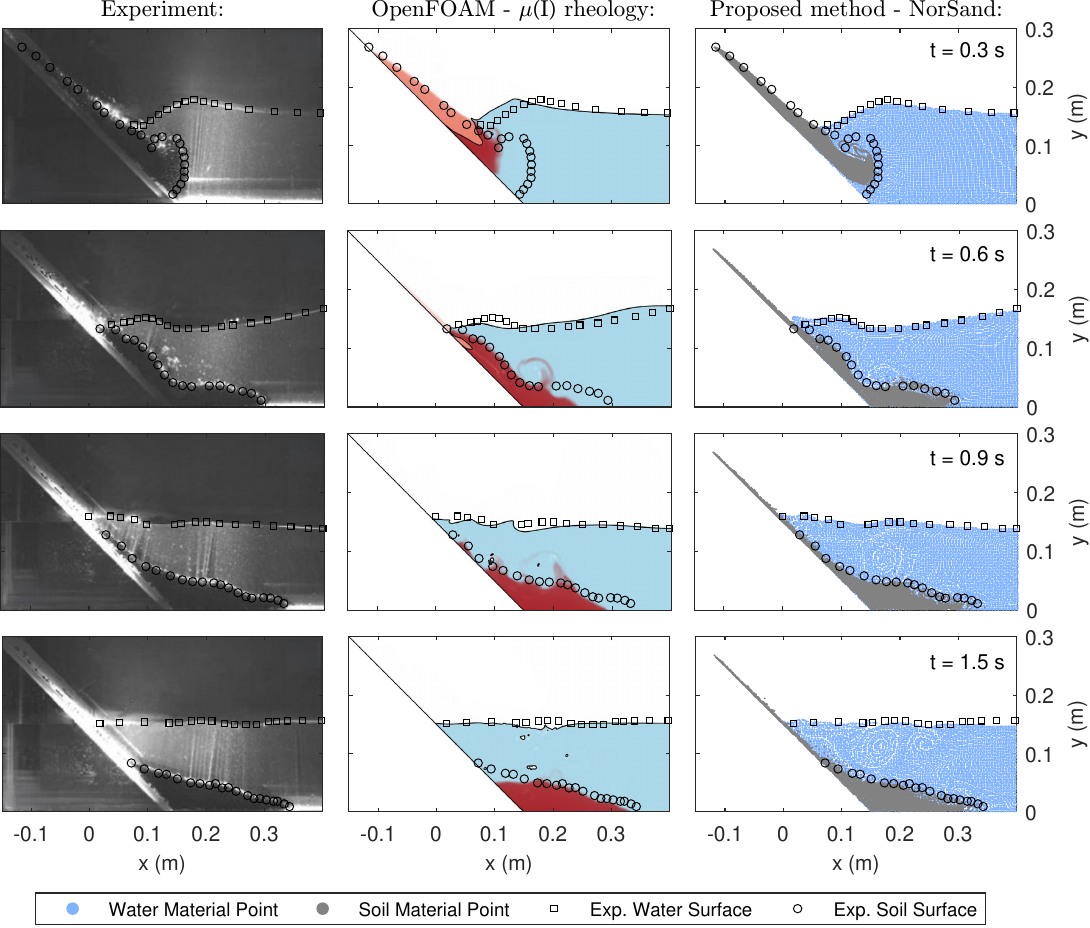}
		\caption{Granular flow and tsunami: case 1 - the evolution of granular flow and generated wave at different time snapshots t = 0.3, 0.6, 0.9, and 1.5~s: comparisons of the
experiment and the OpenFOAM simulation conducted by \citet{Rauter2022} and the proposed method}
		\label{4mmglassbeads_compare}
	\end{center}
\end{figure}

\begin{figure}
	\begin{center}
		\includegraphics[width=0.95\textwidth]{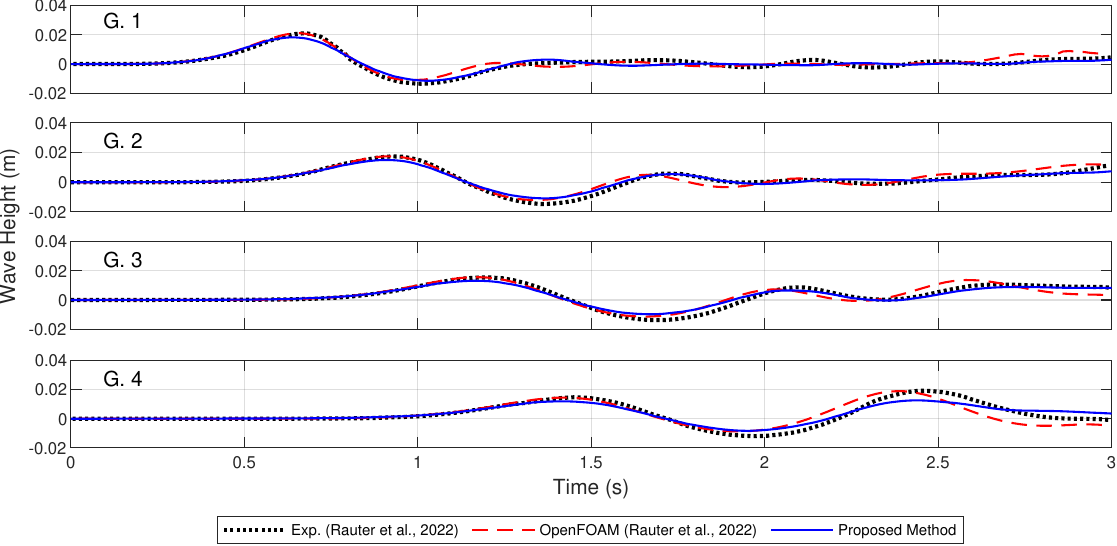}
		\caption{Granular flow and tsunami: case 1 - the wave height time histories recorded at 4 water gauges G. 1, G. 2, G. 3 and G. 4: comparisons of the experiment and the OpenFOAM simulation with $\mu$(I) rheology model conducted by \citet{Rauter2022} and the proposed method with NorSand constitutive model}
		\label{4mmglassbeads_wave}
	\end{center}
\end{figure}

Figure~\ref{4mmglassbeads_compare} compares the simulation given by the proposed numerical model, the experiment and the numerical study conducted by \citet{Rauter2022} using OpenFOAM with $\mu$(I) rheology constitutive model at different time snapshots. The free-water and soil's surfaces from the experiment are plotted as square and circular markers in Fig.~\ref{4mmglassbeads_compare}, respectively. As shown, the proposed method demonstrates superior performance in reproducing the experiment at all time snapshots, compared to the results from OpenFOAM. The NorSand constitutive model can capture the complex geometry and soil behaviour during the landslide process. Due to the accurate prediction of soil behaviour by the proposed numerical method, the water surface also better matches the experiment compared to the previous numerical study conducted by \citet{Rauter2022}, as shown in Fig.~\ref{4mmglassbeads_compare}, despite the previous work having already achieved a good match with the experiment. The excellent agreement with the experiment is also reflected in the wave height time histories recorded at the four water gauges, as shown in Fig.~\ref{4mmglassbeads_wave}. The proposed method also shows better agreement with the experiment compared to that given by OpenFOAM, although the differences are not significant. However, Fig.~\ref{4mmglassbeads_wave} shows that the proposed method very slightly underestimates the wave height compared to the experiment and OpenFOAM. This was also reported in the MPM study using the Drucker-Prager constitutive model conducted by \citet{He2024}. Similar to the $\mu$(I) rheology model, the soil phase governed by the Drucker-Prager constitutive model also showed difficulty matching the experimental results, because neither $\mu$(I) rheology nor the Drucker-Prager model can account for volume changes or density evolution. As the NorSand constitutive model is based on critical state theory, it can reproduce the different responses of soil under both dense and loose conditions.

\clearpage

\subsection{Granular flow and tsunami: case 2}\label{sec:tsunami2}

\begin{figure}[htbp]
	\begin{center}
		\includegraphics[width=0.6\textwidth]{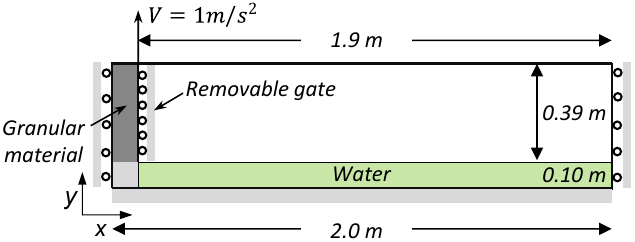}
		\caption{Granular flow and tsunami: case 2 - setup of the initial geometry and boundary conditions}
		\label{5mmglassbeadscollapse_setup}
	\end{center}
\end{figure}

By establishing a physical laboratory model, \citet{Sarlin2022} studied some examples of landslide-induced tsunamis. Figure~\ref{5mmglassbeadscollapse_setup} illustrates one of their models for generating solitary waves. Inside the flume, a granular column (formed of glass beads, 0.39~m in height and 0.1~m in length) is held by a removable gate and a solid base (0.1~m in height and width). The plane strain condition is valid because the flume width (0.1~m) is sufficiently thick. The rest of the flume is filled by a layer of 0.1~m shallow water, which has the same height as the solid base underneath the granular column. Therefore, the glass beads are fully dry at the initial configuration. These glass beads have a 5 mm uniform diameter, and their grain density is 2500~$\text{kg}$/$\text{m}^3$ \citep{Sarlin2022}. According to \citet{Sarlin2022}, at $t = 0$ s the gate was lifted with a 1~m/s constant speed. The authors provided a video with a 50 fps frame rate recording the entire test, from releasing the gate to the granular column collapse, as well as the wave generation due to the impact of the granular flow. The frames at $t = 0.2$, 0.5, 0.8, 1.1 and 1.4~s were extracted from the video and digitised to accurately compare them with the numerical results obtained with the method proposed in this research.

The geometry and boundary conditions for the numerical models follow the setup previously discussed and shown in Fig.~\ref{5mmglassbeadscollapse_setup}. The removable gate is modelled as a roller boundary condition, lifted at a speed of 1~m/s. The initial condition of the model was obtained by gradually applying the gravity load within 1~s before opening the gate. For this numerical example, a structured quadrilateral background mesh with a cell size of 0.01~m is used to discretise the computational domain; $6^2$ and $4^2$ soil and water material points, respectively, are uniformly placed in each cell. Note that a $4^2$ configuration is sufficient for this problem in terms of accuracy. More material points are placed in each cell for a better visualisation purpose. A time increment $\Delta t$ = $6\times10^{-5}$~s (about 0.95 CFL) is used, and the total simulation time is 2~s excluding the time for initialisation.

The simulations were carried out using both NorSand and Mohr-Coulomb constitutive models based on a large-strain framework stabilised by the modified F-bar method developed in \citet{Xie2025}. For the water phase, the TPIC (1\%) scheme with hyperelastic water model and $\delta$-correction method is used. The same NorSand parameters (Table~\ref{tab:4mmglassbeads}) previously calibrated from 4~mm diameter glass beads have been used, assuming similar granular behaviour to that of 5~mm diameter beads, as no triaxial tests for beads of this size were found in the literature. It is well known that the classical Mohr-Coulomb model cannot accurately reproduce the results of triaxial tests. However, it requires only a few parameters, most of which can be obtained directly from the calibrated NorSand model. For example, the shear modulus is set to the same value as NorSand, and the Poisson's ratio is assumed to be 0.3 (typical in granular materials). The friction angle of 25.5° is used in this case. For the granular material, the cohesion is zero. The dilation angle is likewise adopted as zero, resulting in a non-associative flow rule.

\begin{figure}[t]
\centering
\includegraphics[width=0.99\textwidth]{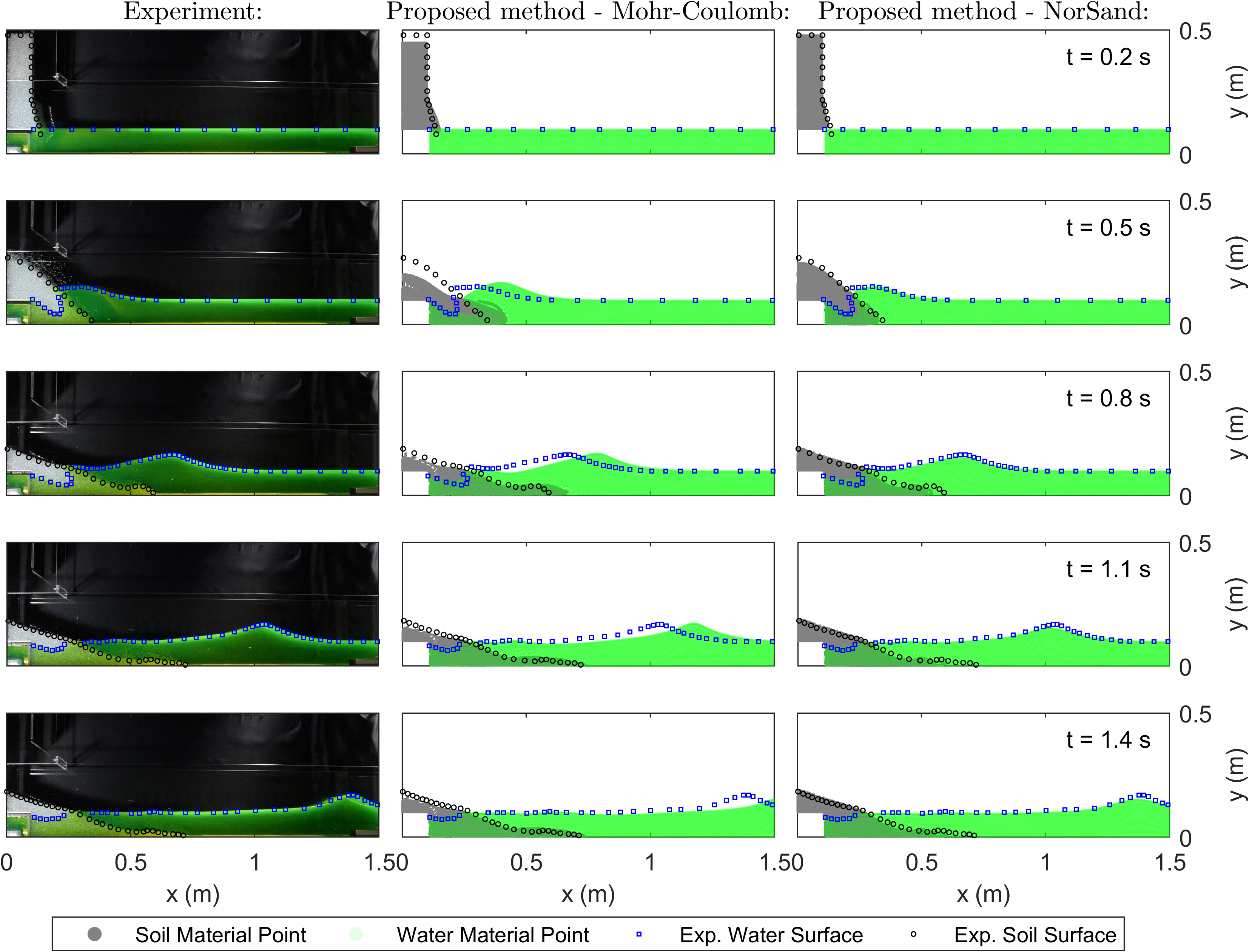}
\caption{Granular flow and tsunami: case 2 - comparisons of the experiment conducted by \citet{Sarlin2022} and the proposed numerical method with Mohr-Coulomb and NorSand constitutive models at different time snapshots.}\label{compare5mmglassbeadsNS}
\end{figure}

Figure~\ref{compare5mmglassbeadsNS} presents the comparisons between the experiment and the simulations obtained with the proposed method for NorSand and Mohr-Coulomb constitutive models. The black circular and blue square markers depict the soil and water surfaces in the experiment images, respectively. As shown in Fig.~\ref{compare5mmglassbeadsNS}, the proposed numerical method with NorSand constitutive model demonstrates a remarkable performance in reproducing the landslide-induced tsunami. Both soil and water phases are very stable under their intense interaction and the process of large deformation. The proposed method with a carefully calibrated NorSand model accurately described the complex geometry and failure process of the granular flow. However, the macroscopic constitutive model cannot describe the rolling of glass bead particles, resulting in an underestimated run-out distance.

Moreover, as shown in Fig.~\ref{compare5mmglassbeadsNS}, the generated tsunami matches perfectly the experiment at all time snapshots, confirming the correctness and stability of the proposed method in simulating wave propagation. After the granular flow impacts the water, it percolates inside the voids of the granular material. Then, the water level within the granular material increased gradually. The proposed semi-implicit double-point MPM with the NorSand model captured this phenomenon accurately, as shown in Fig.~\ref{compare5mmglassbeadsNS}. It is also worth emphasising the excellent stability of the proposed method, as an advanced constitutive model is usually not as stable as classic models like Mohr-Coulomb. The proposed method with the Mohr-Coulomb model fairly matches the experiment, but the agreement is still acceptable, given the low number of parameters needed for this model. In this case, the propagation of the granular flow occurs faster with the Mohr-Coulomb model, and the final run-out distance was also overestimated. Furthermore, the overestimation is evident in the settlement of the crest. The overall behaviour of Mohr-Coulomb is over-softened compared with the experiment and NorSand model. This was expected because this basic model cannot properly represent the dilatancy of the granular material. 

\subsection{Dam break over a movable granular bed}\label{sec:dambreak_movable}

\begin{figure}[htbp]
\centering
\begin{tikzpicture}[scale=6.5]
\draw[thick] (0,0) rectangle (0.90,0.65);

\fill[brown!40] (0,0) rectangle (0.90,0.10);

\fill[blue!25] (0,0.10) rectangle (0.30,0.45);

\draw[dashed, thick] (0.30,0) -- (0.30,0.65);

\draw[<->, >=latex] (0,0.62+0.06) -- (0.90,0.62+0.06);
\node[above] at (0.45,0.62+0.06) {\footnotesize 6.0 m};

\draw[<->, >=latex] (0,0.57) -- (0.30,0.57);
\node[above] at (0.15,0.57) {\footnotesize 3.0 m};

\draw[<->, >=latex] (0.30,0.57) -- (0.90,0.57);
\node[above] at (0.60,0.57) {\footnotesize 3.0 m};

\draw[<->, >=latex] (0.93,0) -- (0.93,0.65);
\node[right] at (0.93,0.325) {\footnotesize 0.50 m};

\draw[<->, >=latex] (-0.03,0.10) -- (-0.03,0.45);
\node[left] at (-0.03,0.275) {\footnotesize $H$ = 0.35 m};

\draw[<->, >=latex] (-0.03,0) -- (-0.03,0.10);
\node[left] at (-0.03,0.05) {\footnotesize 0.10 m};

\node at (0.15,0.28) {\footnotesize \textit{Water}};
\node at (0.60,0.05) {\footnotesize \textit{Saturated sediment bed}};

\node[above left] at (0.88,0.50) {\footnotesize \textit{Computational domain}};

\draw[->, >=latex, thick] (0.02+0.28,0.02-0.02) -- (0.10+0.28,0.02-0.02);
\draw[->, >=latex, thick] (0.02+0.28,0.02-0.02) -- (0.02+0.28,0.10-0.02);
\node[below] at (0.10+0.28,0.02-0.02) {\footnotesize $x$};
\node[left] at (0.02+0.28,0.10-0.02) {\footnotesize $y$};

\fill[pattern=north east lines, pattern color=gray] (0,-0.02) rectangle (0.90,0);
\draw[thick] (0,0) -- (0.90,0);

\foreach \y in {0.05,0.15,0.25,0.35,0.45,0.55}
  \draw (-0.008,\y) circle (0.005);

\foreach \y in {0.05,0.15,0.25,0.35,0.45,0.55}
  \draw (0.908,\y) circle (0.005);

\end{tikzpicture}

\caption{{Dam break over a movable granular bed - setup of the initial geometry and boundary conditions (not to scale)}}
\label{fig:dambreak_movable_setup}
\end{figure}
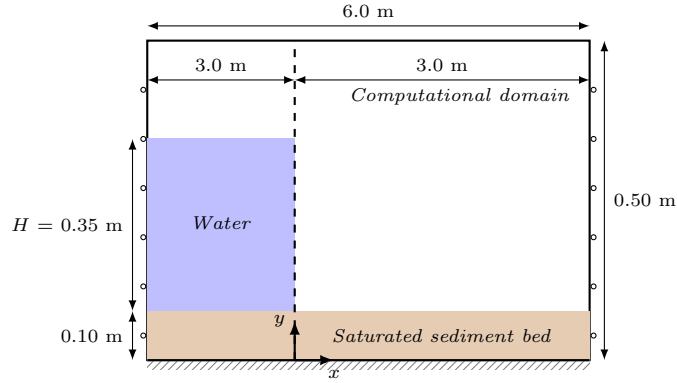

{The final numerical example demonstrates the capability of the proposed method to simulate the reverse interaction: granular deformation driven by a fluid dam break. While the previous examples focused on granular-induced fluid motion (landslide-generated tsunamis), this case examines the erosion and entrainment of bed sediments by a dam-break wave. The proposed numerical method has been compared with the experiment conducted by \citet{Spinewine2005}. In this experiment, a water column of 0.35~m depth is held by a removable gate in a 6~m long flume with a flat bed of loose PVC pellets at the same level on both sides of the gate. After the sudden removal of the gate (within approximately 0.1~s), the dam-break wave propagates downstream, eroding and entraining substantial amounts of bed sediments. The experimental setup is illustrated in Fig.~\ref{fig:dambreak_movable_setup}. \citet{Spinewine2005} used fast digital cameras at 200~fps to record the evolution of the free surface, the bed level and the top of the sediment transport layer through the transparent sidewalls.}

{The PVC pellets used in the experiment have an equivalent spherical diameter of 3.9~mm, an intrinsic density of $\rho_s = 1580$~kg/m$^3$, a friction angle of $\phi_{cs} = 38^{\circ}$, no cohesion, and a random close packing fraction of 0.58. From these properties, the initial porosity is $n_0 = 1 - 0.58 = 0.42$, the initial void ratio is $e_0 = n_0/(1-n_0) = 0.724$, the initial specific volume is $v_0 = 1 + e_0 = 1.724$, and the slope of the critical state line is $M = 6\sin\phi_{cs}/(3-\sin\phi_{cs}) = 1.55$. The initial permeability $k_0$ is estimated from Eq.~\eqref{eq:est_permeability}, using a shape factor $SF = 6.6$ (rounded pellets) and an effective particle diameter $D_{eff} = 0.39$~cm, yielding $k_0 = 0.153$~m/s. Due to the lack of triaxial test data for the 3.9~mm PVC pellets, the NorSand model parameters calibrated for the 4~mm glass beads (Table~\ref{tab:4mmglassbeads}) are adopted, with the exception of the slope of the critical state line $M$, which is updated to 1.55 based on the reported friction angle. The input parameters are summarised in Table~\ref{tab:NorSand_PVC}.}


\begin{table}[htbp]
\centering
{
\caption{Dam break over a movable granular bed - NorSand model parameters for 3.9~mm PVC pellets and input parameters for numerical simulation}
\label{tab:NorSand_PVC}
\begin{tabular}{@{}ll@{}}
\toprule
Parameter & Value \\
\midrule
Water density $\rho_w$ [kg/m$^3$] & 1000 \\
Soil grain density $\rho_s$ [kg/m$^3$] & 1580 \\
Shear modulus $G$ [MPa] & 40 \\
Swelling index $\kappa$ & 0.002 \\
Reference specific volume $v_{c0}$ & 1.84 \\
Compression index $\lambda$ & 0.02 \\
Slope of the critical state line $M$ & 1.55 \\
Yield function constant $N$ & 0.5 \\
Plastic potential constant $\bar{N}$ & 0.5 \\
Hardening coefficient $h$ & 1000 \\
Multiplier for maximum plastic dilatancy $\alpha$ & $-11$ \\
Initial specific volume $v_0$ & 1.724 \\
Initial porosity $n_0$ & 0.42 \\
Initial hydraulic conductivity $k_0$ [m/s] & 0.153 \\
Gravitational acceleration $g$ [m/s$^2$] & $-9.81$ \\
\hline
\end{tabular}
}
\end{table}

{The geometry and boundary conditions for the numerical model closely reproduce the conditions of the experiment. The computational domain has dimensions of 6.0~m in length and 0.5~m in height. A plane strain condition is adopted. Roller and fixed boundary conditions are applied on the sides and bottom of the computational domain, respectively. The sediment bed is 0.10~m thick and extends at the same level on both sides of the gate, with the upstream water depth $H = 0.35$~m above the bed surface. The bed is fully saturated prior to gate removal. Similar to the previous dam-break examples, the gate is not modelled explicitly, but simulated with an immediate release of the constraint of the water in its location; gravitational acceleration $g = 9.81$~m/s$^2$ is applied at $t = 0$~s to generate the dam-break flow. The initial condition of the model was obtained by gradually applying the gravity load before releasing the gate. For this numerical example, a structured quadrilateral background mesh with a cell size of 0.01~m is used to discretise the computational domain; $6^2$ and $4^2$ soil and water material points, respectively, are uniformly placed in each cell at the initial configuration. The TPIC (1\%) scheme with the hyperelastic water model, and the $\delta$-correction method is used for both phases. The modified F-bar method developed in \citet{Xie2025} is used to stabilise the solutions in the soil phase. A time increment $\Delta t = 6 \times 10^{-5}$~s (about 0.95 CFL) is used, and the total simulation time is 1.5~s excluding the time for initialisation.}

\begin{figure}[htbp]
\centering
\includegraphics[width=\textwidth]{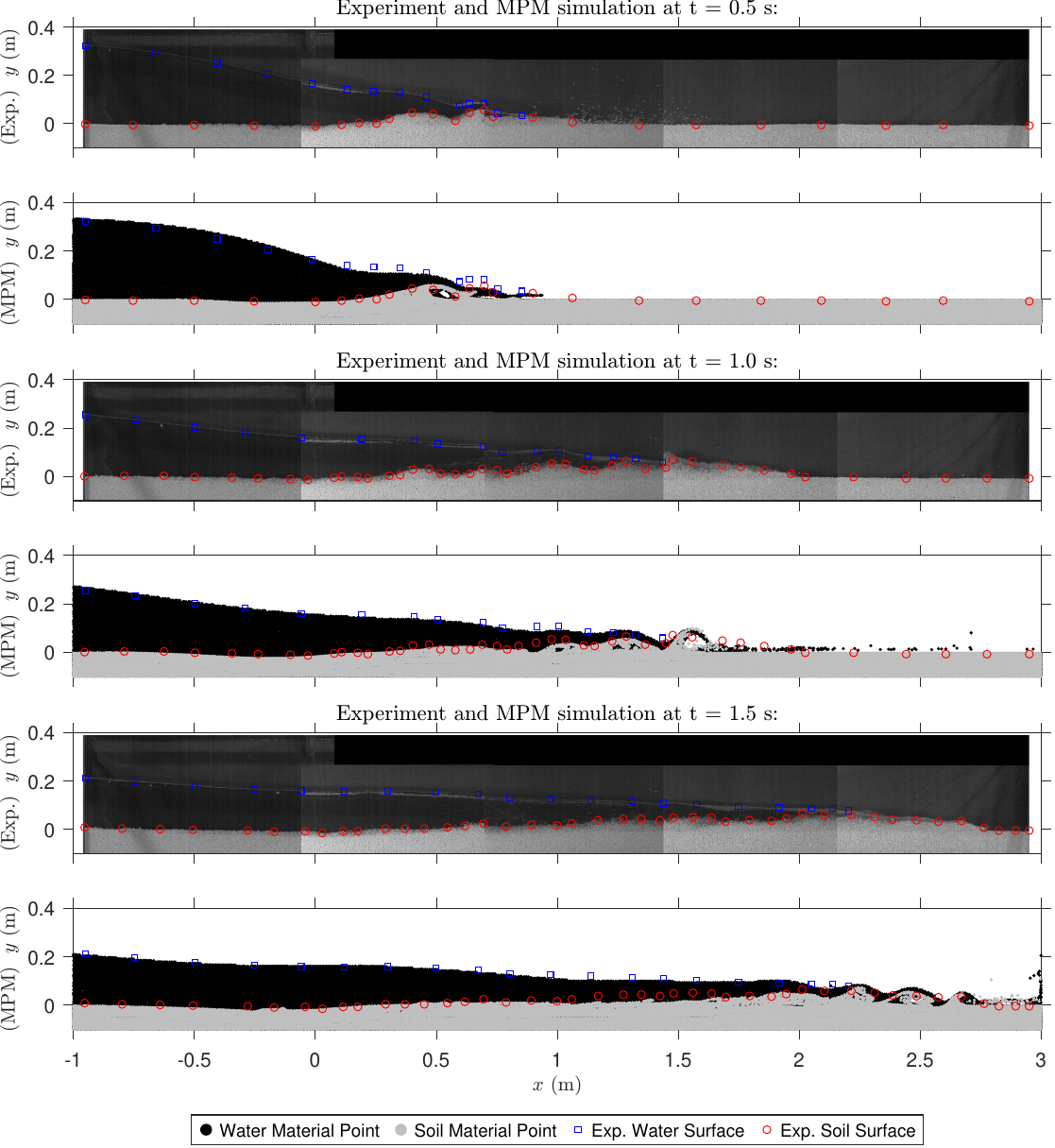}
\caption{{Dam break over a movable granular bed - the 
evolution of the granular flow and free-surface at different time 
snapshots: comparisons of the experiment conducted by 
\citet{Spinewine2005} and the proposed numerical method with NorSand 
constitutive model. Markers denote experimental flow interfaces.}}
\label{fig:dambreak_movable_results}
\end{figure}

{Figure~\ref{fig:dambreak_movable_results} presents the comparison between the experiment and the numerical simulation at different time snapshots. The experimental free surfaces reported by \citet{Spinewine2005} are plotted as markers in the figure to aid their comparison with the numerical results. As shown, the proposed method has a reasonable agreement with the experimental observations. The erosion of the sediment bed by the advancing dam-break wave, the formation of the sediment transport layer, and the evolution of the free surface are all captured by the numerical model. The wavefront celerity and the extent of bed scour near the gate region are also in reasonable correspondence with the experimental data.}

{Some discrepancies between the numerical results and the experiment can be observed, particularly at the wavefront and in the run-out distance of the eroded sediments. These may be attributed to two factors. First, the NorSand model parameters were calibrated for 4~mm glass beads rather than 3.9~mm PVC pellets; dedicated triaxial tests for the PVC material would be needed for a more precise calibration. Second, the macroscopic continuum constitutive model cannot describe the free rolling of individual granular particles, which is a prominent mechanism at the wavefront where the sediment transport layer is fully saturated, and grains move as a dense granular-fluid mixture. Despite these limitations, the results demonstrate that the proposed semi-implicit, double-point MPM can reliably simulate both directions of soil-water interaction: granular-induced fluid motion (as in the tsunami examples) and fluid-induced granular deformation (as in the present case).}

\clearpage

\section{Conclusions}
{This paper proposes several improvements to the semi-implicit, two-phase, double-point MPM for modelling problems involving both free-water and seepage flow. While the overall formulation is presented within the semi-implicit double-point MPM framework, many of the individual contributions have broader applicability beyond this specific method family. In particular, the hyperelastic water model (Sect.~\ref{sec:hyperelastic_water}), the nodal-based free-surface detection method (Sect.~\ref{sec:free_surface_detection}), the smooth porosity-permeability transition via the Kozeny-Carman relation (Sect.~\ref{sec:sw_interface}), and the TPIC(1\%) hybrid velocity integration scheme (Sect.~\ref{sec:1pTPIC}) can be readily adopted in other MPM formulations (single-point, monolithic, explicit, or fully implicit) and, in some cases, in other particle-based methods such as SPH and MPS. The key contributions and findings can be summarised as follows:}
\begin{itemize}
\item A non-linear Darcy law has been implemented to reproduce problems with high Reynolds numbers (fast flows), such as dam break and landslide-induced tsunamis. Despite the higher complexity of this approach, the proposed implementation maintains computational efficiency through the use of velocities from the previous time step in calculating the non-linear term.
\item A hyperelastic constitutive model for water has been successfully integrated into the semi-implicit, coupled MPM for the first time. This model is more accurate than previously reported approaches to represent the slightly compressible nature of water and its viscosity, leading to improved numerical stability, particularly in fast-moving free-water scenarios, as demonstrated in the validation examples and comparisons.
\item A new nodal-based free-water surface detection method has been developed, offering advantages over traditional cell-based approaches. This method allows for stricter tolerance values (0.99 versus 0.75) and produces more accurate and stable results, particularly when using higher-order B-spline shape functions.
\item The implementation of a smooth soil-water interface, combined with the Kozeny-Carman equation for gradual permeability changes, significantly improves the stability of simulations involving soil-water interactions.
\item A hybrid approach, combining TFLIP and TPIC methods (TPIC 1\%), has been proposed and demonstrated to effectively stabilise the velocity fields whilst maintaining energy conservation. This approach, together with the $\delta$-correction method, has proven to be particularly effective in handling complex water-soil interaction problems.
\end{itemize}

{The effectiveness of the combination of these strategies has been assessed through five numerical examples: a classic dam break, a dam break seepage flow, two granular flow-induced tsunamis, and a dam-break wave over a movable granular bed.} The results show excellent agreement with the experimental data reported in the literature, particularly in capturing complex phenomena such as breaking waves and seepage flows. The NorSand constitutive model, when properly calibrated, has demonstrated superior performance in capturing the granular material behaviour compared to simpler models like $\mu$(I) rheology or traditional elasto-plastic models such as  Mohr-Coulomb.

{Compared with other existing MPM formulations for soil-water interaction with free-surface flows, the proposed scheme offers three key advantages: (i) it retains the computational efficiency of the incremental fractional step method, requiring only a scalar pressure Poisson solve per step and, crucially, removing both the permeability-dependent time-step restriction and the volume conservation issue inherent in traditional fractional step approaches, which otherwise require very small time steps for low-permeability geotechnical problems and experience progressive loss of water volume over long simulations; (ii) the combination of stabilisation ingredients (hyperelastic water model, TPIC(1\%), nodal free-surface detection, smooth interface, and $\delta$-correction) enables stable simulation of problems involving simultaneous free-surface flow, seepage, and intense soil-water interaction, which remain challenging for existing semi-implicit MPM formulations; and (iii) the use of an advanced critical-state constitutive model (NorSand), rather than simplified elasto-plastic models employed in previous double-point MPM studies, provides a more physically meaningful description of granular materials undergoing large deformation and phase interaction. The final numerical example further confirms that the proposed formulation captures both directions of soil-water interaction: granular-induced fluid motion and fluid-induced granular deformation.}

These developments represent a significant step forward in the capability of MPM to handle complex geotechnical problems involving free water and soil-water coupling. {Future work will focus on developing a fully implicit, double-point MPM for soil-water coupled problems with free-surface flows, which, to the best of our knowledge, does not yet exist in the literature. The proposed formulation can also be applied to a broader 
range of practical geotechnical engineering problems, including real-scale landslide-generated tsunamis, breach of earth dams, and other cascading geohazards.} Finally, the robust numerical formulation and sufficient validations presented in this paper pave the way to simulate real landslides and cascading disasters. In such cases, the model parameters will have to be accurately calibrated using laboratory tests to evaluate the stress-strain behaviour for a wide range of strains until the critical state is reached. Ideally, triaxial tests would be performed on soil samples collected on-site.

\backmatter








\bmhead{CRediT author statement}
\textbf{Mian Xie}: Conceptualisation, Methodology, Software, Validation, Visualisation, Writing - Original Draft. \textbf{Pedro Navas}: Conceptualisation, Supervision, Writing - Review \& Editing. \textbf{Susana L\'opez-Querol}: Conceptualisation, Supervision, Writing - Review \& Editing.

\section*{Declarations}
\bmhead{Conflict of interest}
The authors declare no potential conflict of interest.

\bmhead{Data availability}
Data are available by the corresponding author after reasonable request.

\begin{appendices}






\end{appendices}


\bibliography{sn-bibliography}

\end{document}